\documentclass[reqno,11pt]{amsart}
\usepackage{amsthm}
\usepackage{amsmath,amsfonts,amssymb,bm,bbm}
\usepackage{xr}
\usepackage{mathrsfs}
\usepackage{mathtools}
\usepackage{hyperref}
\usepackage{mathabx}
\usepackage[initials,nobysame]{amsrefs} 

\numberwithin{equation}{section}

\hypersetup{hypertex=true,
colorlinks=true,
linkcolor=blue,
anchorcolor=blue,
citecolor=blue}
\newtheorem{theorem}{Theorem}[section]
\newtheorem{lemma}[theorem]{Lemma}
\newtheorem{proposition}[theorem]{Proposition}

\newtheorem{remark}[theorem]{Remark}
\usepackage[top=2.0cm,bottom=2.0cm,left=3cm,right=3cm]{geometry}
\numberwithin{equation}{section}
\mathtoolsset{showonlyrefs}
\begin{document}
\title[Kinetic-fluid boundary layers and acoustic limit]{Kinetic-fluid boundary layers and acoustic limit for the Boltzmann equation with general Maxwell reflection boundary condition}
\date{}
\author[N. Jiang and Y.-L. Wu]{Ning Jiang \and Yulong Wu}

\address[N. Jiang]{School of Mathematics and Statistics, Wuhan University, Wuhan 430072, China
}

\email{
njiang@whu.edu.cn}
\address[Y.-L. Wu]{
  School of Mathematics and Statistics,
    Wuhan University,  Wuhan 430072, China
}
  \email{
  yulong\string_wu@whu.edu.cn
}
\maketitle

\begin{abstract}
 We prove the acoustic limit from the Boltzmann equation with hard sphere collisions and the Maxwell reflection boundary condition. Our construction of solutions include the interior fluid part and Knudsen-viscous coupled boundary layers. The main novelty is that the accommodation coefficient is in the full range $0<\alpha\leq 1$. The previous works in the context of classical solutions only considered the simplest specular reflection boundary condition, i.e. $\alpha=0$. The mechanism of the derivation of fluid boundary conditions in the case $\alpha=O(1)$ is quite different with the cases $\alpha=0$ or $\alpha=o(1)$.   This rigorously justifies the corresponding formal analysis in Sone's books \cite{sone2002kinetic,sone2007molecular}. In particular, this is a smooth solution analogue  of \cite{jiang2010remarks}, in which the renormalized solution was considered and the boundary layers were not visible.
  \end{abstract}

%\tableofcontents(生成目录)

\section{Introduction and main result}
\subsection{Problem presentation}
In this paper, we study the acoustic limit for the Boltzmann equation in half-space with the Maxwell reflection boundary condition. The following Boltzmann equation with Euler scaling is considered:
\numberwithin{equation}{section}%%%%%%%%%%%%%%%%%%%%%%%%%%%%%%%%%%%编号有章节
\begin{equation}
\begin{cases}\label{1.1.1}
\partial_tF_\varepsilon+v\cdot\nabla_xF_\varepsilon=\frac 1\varepsilon Q(F_\varepsilon,F_\varepsilon)  &\text{on} \quad \mathbb{R}_+\times\mathbb{R}^3_+\times\mathbb{R}^3, \\
F_\varepsilon(0,x,v)=F^{in}_\varepsilon(x,v)\geq 0     &\text{on}\quad \mathbb{R}^3_+\times\mathbb{R}^3.\\
\end{cases}
\end{equation}
Here, $F_\varepsilon=F_\varepsilon(t,x,v) \geq0$ denotes the density distribution function of the particle gas at time $t \geq0$, position $x\in \mathbb{R}^3_+$, and velocity $v\in \mathbb{R}^3$. The positive parameter $\varepsilon$ is the Knudsen number defined as the ratio of the molecular mean free path length to a representative physical length scale. Moreover, let $v,v_\ast $ and $v',v'_\ast $ be the velocities of the particles before and after the collision respectively, which satisfy
 $$
\begin{aligned}
&\mathrm{a}.\, \text{Conservation laws } v'=v+[(v_\ast -v)\cdot w]w,\quad|v_\ast |^2+|v|^2=|v'_\ast |^2+|v'|^2,\\
&\mathrm{b}.\, \text{The formulation } v_\ast'=v_\ast-[(v_\ast -v)\cdot w]w \quad  \text{ for } \omega \in \mathbb{S}^2. 
\end{aligned}
$$
Throughout this paper, the Boltzmann collision operator $Q(F_1,F_2)$ takes the form
\begin{equation}\label{1.1.2}
Q(F_1,F_2)=\int_{\mathbb{R}^3 \times \mathbb{S}^2}b(|v_\ast-v|,\omega )[F_1(v_\ast^{\prime})F_2(v^{\prime})-F_1(v_\ast)F_2(v)]dv_\ast d\omega,
\end{equation}
where
$$b(|v_\ast-v|,\omega )=|(v_\ast-v)\cdot\omega|.$$

Let $\Sigma:=\partial \mathbb{R}^3_+ \times \mathbb{R}^3$ be the phase space boundary of $\mathbb{R}^3_+ \times \mathbb{R}^3$. The phase boundary $\Sigma$ can be split by outgoing boundary $\Sigma_+$, incoming boundary $\Sigma_-$, and grazing boundary $\Sigma_0$:
$$\Sigma_+=\{(x,v):x \in \partial\mathbb{R}^3_+,v\cdot n>0\},$$
$$\Sigma_-=\{(x,v):x \in \partial\mathbb{R}^3_+,v\cdot n<0\},$$
$$\Sigma_0=\{(x,v):x \in \partial\mathbb{R}^3_+,v\cdot n=0\},$$
where $n$ is the unit outer normal vector $(0,0,-1)$. Let $\gamma_{\pm}F=\mathbbm{1}_{\Sigma_{\pm}}F$.
Moreover, the equation \eqref{1.1.1} is imposed with the Maxwell reflection boundary condition
\begin{equation}\label{1.1.3}
\gamma_-F_\varepsilon=(1-\alpha)L\gamma_+F_\varepsilon+\alpha K\gamma_+F_\varepsilon
 \quad \text{on} \quad \mathbb{R}^+\times \Sigma_-.
\end{equation}
The specular-reflection $L\gamma_+F_\varepsilon$ and the diffuse-reflection part $K\gamma_+F_\varepsilon$ in \eqref{1.1.3} are
\begin{equation}
  L\gamma_+F_\varepsilon(t,x,v)=F_\varepsilon(t,x,R_x v),
\end{equation}
\begin{equation}
  K\gamma_+F_\varepsilon(t,x,v)=M^w\int_{v'\cdot n>0}{(v'\cdot n) \gamma_+F_\varepsilon} d v',
\end{equation}
respectively, where $R_x v=v-2[(v-u_w)\cdot n]n$ and $M^w$ defined by
\begin{equation}
  M^w=\frac {1}{(2\pi {\theta_w})^{3/2}}\exp({-\frac{|v-u_w |^2}{2\theta_w }})
\end{equation}
is a Maxwellian distribution having the wall temperature $\theta_w$ and wall velocity $u_w$. Throughout this paper we assume $u_w=0$ and $\theta_w=1$ at the wall. The nonnegative constant $\alpha$ is the so-called accommodation coefficient, which describes how much the molecules accommodate to the state of the wall. Here we assume that $0<\alpha\leq 1$ and $\alpha=O(1)$, so that the case of specular or almost specular reflection is excluded and diffusive reflection is included. It's worth to emphasize that the Maxwell reflection boundary condition \eqref{1.1.3} ensures no instantaneous mass flux across the boundary:
\begin{equation}
  \int_{\mathbb{R}^3}(v\cdot n)F_\varepsilon(t,x,v)dv=0\quad \text{for all } x\in \partial\mathbb{R}^3_+.
\end{equation}

\subsection{History and motivation}

There have been a lot of significant mathematical studies on the hydrodynamic limit from the Boltzmann equation to the fluid equations, such as incompressible and compressible Navier-Stokes, Euler equations. Starting from the late 1980s, the so-called BGL program aimed at justifying the limit to Leray's solutions of the incompressible Navier-Stokes equations from the Boltzmann equation was initialized by Bardos-Golse-Levermore \cite{bardos1993fluid,bardos1991fluid} and completed by Golse and Sait-Raymond \cite{golse2004navier,golse2009incompressible}. The corresponding results in domain with boundary were carried out in \cite{jiang2017boundary,masmoudi2003boltzmann}. For incompressible Euler limits, we also mention the work of Saint-Raymond \cite{saint2003convergence,saint2009hydrodynamic}. Another group of results is in the framework of classical solutions, which are based on nonlinear energy method, semi-group method or hypercoercivity (the
latter two further rely on the spectral analysis of the linearized Boltzmann operator), see \cite{bardos1991classical,briant2015boltzmann,briant2019boltzmann,gallagher2020convergence,guo2006boltzmann,jang2021incompressible,jiang2018incompressible}.

In the domain with boundary, the situation is much more complicated. Japanese school made significant contributions on formal analysis and numerical studies on the hydrodynamic limits from the Boltzmann equation with Maxwell reflection boundary condition and incoming boundary condition by using Hilbert-type expansions (see \cite{sone2002kinetic,sone2007molecular,aoki2017slip,hattori2019sound,takata2012asymptotic}).  We would like to mention that the mathematical analysis of viscous and kinetic layers has already been used in the work of the first author of the present paper and Masmoudi \cite{jiang2017boundary}, in which the viscous and kinetic layers are used to deal with the acoustic waves. Later, Guo-Huang-Wang \cite{guo2021hilbert} started to rigorously justify the compressible Euler and acoustic limits using the formal derivation in Sone's books \cite{sone2007molecular,sone2002kinetic}, in which the specular reflection boundary condition, i.e., the accommodation coefficient $\alpha=0$ was considered. Furthermore, the first author of the current paper, Luo and Tang \cite{jiang2021compressiblea} considered the case $\alpha=O(\sqrt{\varepsilon})$. Both fluid and kinetic layers are needed in the cases $\alpha=0$ and $\alpha=O(\sqrt{\varepsilon})$. We point out that all the boundary layers in \cite{guo2021hilbert,jiang2021compressiblea,jiang2021compressibleb} are linear because they appear in higher order terms. Specifically, if $\alpha=O(1)$ the leading term of the fluid layer satisfies the nonlinear compressible Prandtl equations, whose well-posedness (even local) in Sobolev-type space is completely open. It should be a challenging and interesting problem to consider the compressible Euler limit from the Boltzmann with $\alpha=O(1)$. As for the acoustic limit, \cite{golse2002stokes} and \cite{jiang2010remarks} worked on the acoustic limit in the periodic spatial domain and bounded domain with  Maxwell reflection boundary condition, respectively, in the framework work of renormalized solutions. Later, Guo-Jang and the first author of this paper \cite{guo2010acoustic} worked on the classical solution in the domain without boundary, in which the Hilbert expansion method and $L^2-L^\infty$ framework was employed.

In the current paper, we treat the acoustic limit with the accommodation coefficient $\alpha=O(1)$ with the expansion around the global Maxwellian $\mu$. To the best of our knowledge, there is no rigorously mathematical analysis in the hydrodynamic limit of the Boltzmann equation with the accommodation coefficient $\alpha=O(1)$. One of the main reasons is that the derivation of the boundary condition of the fluid equations is quite different with the cases of specular reflection \cite{guo2021hilbert} with $\alpha=0$ and almost specular reflection \cite{jiang2021compressiblea}.  In these cases, the number of the solvability condition of the Knudsen-layer problem is four (see Golse-Perthame-Sulem \cite{golse1988boundary}). This provides four boundary conditions to the corresponding fluid equations. However, the compressible Euler or acoustic system requires only one boundary condition. This mismatch indicates the existence of the viscous layer, i.e. Prandtl layer. For the case $\alpha=O(1)$, the mechanism of the derivation of fluid boundary conditions are quite different. When $\alpha=O(1)$, the number of the solvability condition of the Knudsen layer problem is only one. This seems to indicate that the viscous layer is not needed. However, there is an extra condition: the solutions of the Knudsen layer equation have to vanish at infinity. We must utilize the symmetric properties of the Knudsen layer problem based on the properties of the linearized Boltzmann operator $\mathcal{L}$ and the boundary operator (see \cite{aoki2017slip,hattori2019sound,takata2012asymptotic,sone2002kinetic,sone2007molecular} for the formal analysis and \cite{Jiang2024boundary} for the rigorous proof). We still need four boundary conditions for the fluid equations to make the solution of the Knudsen layer vanishes at infinity.

We emphasize that this approach is not new. Sone and his co-workers have already used the symmetrics properties to investigate the Knudsen layer solution in the Hilbert expansion method and derived four boundary conditions to the fluid equations, see \cite{sone2002kinetic,sone2007molecular,aoki2017slip,hattori2019sound,takata2012asymptotic} and references therein. Recently, the authors of this paper and He proved this statement rigorously, see \cite{Jiang2024boundary}. In summary, the Knudsen layer problem still provides four boundary conditions to the corresponding fluid system in order to vanish at infinity. However, the hyperbolic system \eqref{1.4.1.12} only needs one boundary condition. Thereby the viscous boundary layer is naturally needed. We finally derive the Dirichlet boundary condition for the viscous boundary layer problem, while the corresponding boundary condition is the Neumann boundary condition in \cite{guo2021hilbert} and the Robin boundary condition in \cite{jiang2021compressiblea}. More details will be given in subsection \ref{se1.4.3}. This is the first rigorous justification result in the compressible fluid models with boundary from the Boltzmann equations for the accommodation coefficient $\alpha=O(1)$, although only for the linear acoustic system.

\subsection{Notations}
Throughout this paper, we denote $\bar{U}=(U_1,U_2)$ for any vector $U=(U_1,U_2,\allowbreak U_3)\in \mathbb{R}^3$. Moreover, the notation $V^0:=V|_{x_3=0}$ indicates the value on the boundary for any function $V=V(x)$.

For functional spaces, $H^s$ denotes the standard Sobolev space $W^{s,2}(\mathbb{R}^3_+)$ with norm $\|\cdot\|_{H^s},$ $\|\cdot\|_2$ and $\|\cdot\|_\infty$ denote the $L^2$-norm and $L^\infty$-norm over $(x,v)\in\mathbb{R}^3_+\times\mathbb{R}^3$ variables, respectively. Moreover, $\left<\cdot,\cdot\right>$ denotes the $L^2$-inner product over $(x,v)\in\mathbb{R}^3_+\times\mathbb{R}^3$.

For any function $G=G(t,\bar{x},y,v)$, with $(t,\bar{x},y,v)\in \mathbb{R}_+\times\mathbb{R}^2\times\bar{\mathbb{R}}^3_+\times \mathbb{R}^3$, the Taylor expansion at $y=0$ is
\begin{equation*}
  G=G^0+\sum_{1\leq l\leq N}\frac{y^l}{l!}G^{(l)}+\frac{y^{N+1}}{(N+1)!}\tilde{G}^{(N+1)}.
\end{equation*}
Here the symbols
\begin{equation*}
  G^0=G|_{x_3=0},\quad   G^{(l)}=(\partial^l_yG)(t,\bar{x},0,v),
\end{equation*}
and 
\begin{equation*}
  \tilde{G}^{(N+1)}=(\partial^{N+1}_yG)(t,\bar{x},\eta,v) \quad \text{ for some }\eta\in(0,y).
\end{equation*}
For later use, we define the following linearized Boltzmann operator $\mathcal{L}$:
\begin{equation}\label{1.3.1}
   \mathcal{L} f(v)=-\frac{1}{\sqrt{\mu}}\{Q(\mu,\sqrt{\mu }f)+Q(\sqrt{\mu}f,f)\},
  \end{equation}
  where the global Maxwellian $\mu$ is defined by

\begin{equation}\label{1.3.2}
   \mu(v)=\frac{1}{(2\pi)^\frac{3}{2}}e^{-\frac{|v|^2}{2}}.
\end{equation}
  The bilinear term is defined by
  \begin{equation}\label{1.3.3}
    \Gamma(f,g)=\frac{1}{\sqrt{\mu}}Q(\sqrt{\mu}f,\sqrt{\mu}g).
  \end{equation}
It's well known that the null space $\mathcal{N}$ of $\mathcal{L}$ is generated by (see \cite{caflisch1980fluid} for instance)
\begin{equation}\label{1.3.4}
  \begin{aligned}
   & \chi_0(v)=\sqrt{\mu},\\
    &\chi_i(v)=v_i\sqrt{\mu} \quad (i=1,2,3),\\
   & \chi_4(v)=\frac{|v|^2-3}{2}\sqrt{\mu}.\\
  \end{aligned}
\end{equation}
Clearly, $ \mathcal{L}$ can be decomposed as $$ \mathcal{L}=\nu(v)-K,$$
where $K$ is a self-adjoint operator, given by
\begin{equation}\label{1.3.7}
Kf=\int_{\mathbb{R}^3}k(v,v')f(v')dv',
\end{equation}
and $\nu(v)$ denotes the collision frequency, defined as
\begin{equation*}
  \nu(v)=\int_{\mathbb{R}^3}\int_{\mathbb{S}^2}|(v-v_*)\cdot \omega|\mu(v)d\omega dv_*,
\end{equation*}
with bounds given by
\begin{equation}\label{1.3.5}
  \nu_0(1+|v| ) \le \nu(v) \le \nu_1(1+|v|),
\end{equation}
where $\nu_0,\nu_1$ are given positive constants. Furthermore, we introduce the $L^2_v$ projection $\mathbf{P}$ with respect to $\mathcal{L}$. Then there exists a positive constant $c_0>0$ such that
\begin{equation}\label{1.3.6}
  \left<\mathcal{L}g,g\right>\geq c_0\|(\mathbf{I-P})g\|_\nu^2,
\end{equation}
where the weight $L^2$-norm is
\begin{equation*}
  \|g\|_\nu^2:=\int_{\mathbb{R}^3_+\times\mathbb{R}^3} g^2(x,v)\nu dxdv.
\end{equation*}

%Properties on $\Gamma(f,f)$:
%\begin{equation}
 % \Gamma(f,f) \in \mathcal{N}^{\perp} \text{ and }\left\|\nu^{-1}w\Gamma(f,f)\right\|_{L^\infty_v}\leq C\left\|wf\right\|^2_{L^\infty_v}.
%\end{equation}

In order to quantitatively describe the linear hyperbolic system \eqref{1.4.1.12} in the half-space $\mathbb{R}_+^3$, we let
$$\partial^\alpha_{t,\bar{x}}=\partial_t^{\alpha_0}\partial_{x_1}^{\alpha_1}\partial_{x_2}^{\alpha_2},$$
where $\alpha=(\alpha_0,\alpha_1,\alpha_2)\in \mathbb{N}^3$, and let
\begin{equation}\label{1.3.10}
 \|f(t)\|_{\mathcal{H}^k(\mathbb{R}^3_+)}=\sum_{|\alpha|+i\leq k}\|\partial^\alpha_{t,\bar{x}}\partial^i_{x_3}f(t)\|_{L^2(\mathbb{R}^3)},\quad \|g(t)\|^2_{\mathcal{H}^k(\mathbb{R}^2)}=\sum_{|\alpha|\leq k}\|\partial^\alpha_{t,\bar{x}}g(t)\|^2_{L^2({\mathbb{R}^2})}
\end{equation}
for functions $f(t)=f(t,\bar{x},x_3)$ and $g(t)=g(t,\bar{x})$. We note that $\|f\|_{\mathcal{H}^k}$ is equivalent to the standard Sobolev norm $\|f\|_{\mathcal{H}^k(\mathbb{R}^3)}$ for a function $f$ independent of the variable $t$.

While characterizing quantitatively the heat equations appear in the viscous boundary layer, some new norms are required to be introduced. For $l\geq 0$, we define the weighted norm as follows:
\begin{equation}\label{1.3.11}
 \|f\|^2_{L^2_l}=\int_{\mathbb{R}^2}\int_{\mathbb{R}_+}(1+\zeta)^l|f(\bar{x},\zeta)|^2d\bar{x}d\zeta.
\end{equation}
We further introduce a weighted Sobolev space $\mathbb{H}^r_l(\mathbb{R}^3_+)$ for any $r,l\geq 0$. We denote the multi-index $\beta=(\beta_1,\beta_2)\in \mathbb{N}^2$. For any $l,r\geq 0$, let
\begin{equation}\label{1.3.12}
 l_j=l+2(r-j), 0\leq j\leq r.
\end{equation}
We introduce the norms
\begin{equation}\label{1.3.13}
 \begin{array}{l}
   \sum_{j=0}^{r}\|f(t)\|_{l, j, n}^{2}=
   \sum_{j=0}^{r}\sum_{2 \gamma+|\beta|=j-n}\left\|\partial_{t}^{\gamma} \partial_{\bar{x}}^{\beta} \partial_{\zeta}^{n} f(t)\right\|_{L_{l_{j}}^{2}}^{2}\ \ (0 \leq n \leq j), \\
   \|f(t)\|_{l, r}^{2}=\sum_{n=0}^{r}\|f(t)\|_{l, r, n}^{2}=\sum_{2 \gamma+|\beta|+n=r}\left\|\partial_{t}^{\gamma} \partial_{\bar{x}}^{\beta} \partial_{\zeta}^{n} f(t)\right\|_{L_{l_{r}}^{2}}^{2}, \\
   \|f(t)\|_{\mathbb{H}_{l, n}^{r}\left(\mathbb{R}_{+}^{3}\right)}^{2}=\sum_{j=0}^{r}\|f(t)\|_{l, j, n}^{2}\ \ (0 \leq n \leq j \leq r), \\
   \|f(t)\|_{\mathbb{H}_{l}^{r}\left(\mathbb{R}_{+}^{3}\right)}^{2}=\sum_{j=0}^{r}\|f(t)\|_{l, j}^{2}=\sum_{n=0}^{r}\|f(t)\|_{\mathbb{H}_{l, n}^{r}\left(\mathbb{R}_{+}^{3}\right)}^{2}
   \end{array}
\end{equation}
for function  $f=f(t, \bar{x}, \zeta)$ . For  $g=g(\bar{x}, \zeta)$ , let
\begin{equation}\label{1.3.14}
 \|g\|_{\mathbb{H}_{l}^{r}\left(\mathbb{R}_{+}^{3}\right)}^{2}=\sum_{j=0}^{r} \sum_{|\beta|+n=j}\left\|\partial_{\bar{x}}^{\beta} \partial_{\zeta}^{n} g\right\|_{L_{l_{j}}^{2}}^{2} .
\end{equation}
Similarly, for function $ h=h(t, \bar{x})$ , let
\begin{equation}\label{1.3.15}
 \|h(t)\|_{\mathbb{H}^{r}\left(\mathbb{R}^{2}\right)}^{2}=\sum_{j=0}^{r} \sum_{2 \gamma+|\beta|=j}\left\|\partial_{t}^{\gamma} \partial_{\bar{x}}^{\beta} h(t)\right\|_{L^{2}\left(\mathbb{R}^{2}\right)}^{2} .
\end{equation}

\subsection{Main results}
Motivated by Caflish's work \cite{caflisch1980fluid}, in the present paper we aim to rigorously justify the acoustic limit from the scaled Boltzmann equation by the Hilbert expansion approach. As shown in subsection \ref{subse:truncate}, we truncate the Hilbert expansion as follows:
  \begin{equation*}
        F_\varepsilon=\mu+\sum_{k=1}^{6}\sqrt{\varepsilon}^kF_k(t,x,v)+\sum_{k=1}^{6}\sqrt{\varepsilon}^kF^b_k(t,\bar{x},\zeta,v)
         +\sum_{k=2}^{6}\sqrt{\varepsilon}^kF^{bb}_k(t,\bar{x},\xi,v)+\sqrt{\varepsilon}^5\mu^{\frac{1}{2}}f_{R,\varepsilon}\geq 0,
    \end{equation*}
where the global Maxwellian $\mu$ is defined in \eqref{1.3.2}.

 We now clarify the initial data of the scaled Boltzmann equation \eqref{1.1.1}. For $1\leq k\leq 6$, $F^{in}_k(x,v),F^{b,in}_k(\bar{x},\zeta,v),F^{bb,in}_k(\bar{x},\xi,v)$ can be constructed by the same ways of constructing the expansions  $F_k(t,x,v),F^{b}_k(t,\bar{x},\zeta,v),F^{bb}_k(t,\bar{x},\xi,v)$ in section \ref{Se1.4}, respectively. We impose the well-prepared initial data on the scaled Boltzmann equation
 \begin{equation}\label{1.5.11}
  \begin{aligned}
      F_\varepsilon(0,x,v)=\mu+\sum_{k=1}^{6}\sqrt{\varepsilon}^kF^{in}_k(x,v)+\sum_{k=1}^{6}\sqrt{\varepsilon}^kF^{b,in}_k(\bar{x},\zeta,v)\\
       +\sum_{k=2}^{6}\sqrt{\varepsilon}^kF^{bb,in}_k(\bar{x},\xi,v)+\sqrt{\varepsilon}^5\mu^{\frac{1}{2}}f^{in}_{R,\varepsilon}(x,v).
  \end{aligned}
  \end{equation}
We further define
\begin{equation}\label{1.5.12}
  h_{R,\varepsilon}=w_l(v)f_{R,\varepsilon},
\end{equation}
with the velocity weighted function
\begin{equation}\label{1.5.13}
  w_l(v):=\{1+|v|^2\}^{\frac{l}{2}} \quad \text{for} \quad l\geq0.
\end{equation}

Now, we state our main results on the acoustic limit from the Boltzmann equation.
\begin{theorem}\label{Thm1.4}
  Recall the definitions of $\mathcal{H}^k(\mathbb{R}^3_+)$ and $\mathbb{H}^r_l(\mathbb{R}^3_+)$ in \eqref{1.3.10} and \eqref{1.3.14}. Let $w_l(v)$ be defined in \eqref{1.5.13} with $l\geq 7$ and integers $s_k,s^b_k,s^{bb}_k,l^b_k(1\leq k\leq 6)$ be described as in Proposition \ref{Pro2.4}. Assume that
 \begin{equation}\label{1.5.14}
  \mathcal{E}^{in}:=\sum_{k=1}^{6}\left\{\|(\rho^{in}_k,u^{in}_k,\theta^{in}_k)\|_{\mathcal{H}^{s_k}(\mathbb{R}^3_+)}+\|(\bar{u}^{b,in}_k,\theta^{b,in}_k)\|_{\mathbb{H}^{s^b_k}_{l^b_k}(\mathbb{R}^3_+)}\right\} < \infty.
 \end{equation}
 For $1\leq k\leq 6$, let $F_k(t,x,v),F^b_k{(t,\bar{x},\zeta,v)},F^{bb}_k{(t,\bar{x},\xi,v)}$ be constructed in  section \ref{Se1.4}. Suppose that
 \begin{equation}
  \mathcal{E}^{in}_R:=\sup_{\varepsilon\in(0,\varepsilon_0)}\left\{\|f_{R,\varepsilon}^{in}\|_2+\sqrt{\varepsilon}^3\|h^{in}_{R,\varepsilon}\|_\infty\right\} < \infty.
 \end{equation}
 Let $\tau>0$ be any prescribed constant. There is a small constant $\varepsilon_0>0$ such that the scaled Boltzmann equation \eqref{1.1.1} with Maxwell reflection boundary condition \eqref{1.1.3} and well-prepared initial data \eqref{1.5.11} admits a unique solution for $\varepsilon\in(0,\varepsilon_0)$ over the time interval $t\in[0,\tau]$ with the expanded form \eqref{1.5.1}, i.e.,

  \begin{equation}
    \begin{aligned}
        F_\varepsilon=\mu+\sum_{k=1}^{6}\sqrt{\varepsilon}^kF_k(t,x,v)+\sum_{k=1}^{6}\sqrt{\varepsilon}^kF^b_k(t,\bar{x},\zeta,v)\\
         +\sum_{k=2}^{6}\sqrt{\varepsilon}^kF^{bb}_k(t,\bar{x},\xi,v_3)+\sqrt{\varepsilon}^5\mu^{\frac{1}{2}}f_{R,\varepsilon}\geq 0,
    \end{aligned}
    \end{equation}
where the remainder $f_{R,\varepsilon}$ satisfies
\begin{equation}
  \sup_{t\in [0,\tau]}\left(\|f_{R,\varepsilon}\|_2+\sqrt{\varepsilon}^3\|h_{R,\varepsilon}\|_\infty\right)\leq C(\tau,\mathcal{E}^{in},\mathcal{E}^{in}_R)<\infty.
\end{equation}
\end{theorem}

\begin{remark}
  Combining with the Proposition \ref{Pro2.4}, Theorem \ref{Thm1.4} indicates that as $\varepsilon\to 0$,
  \begin{equation}
    \sup_{t\in [0,\tau]}\left(\left\|\frac{F_\varepsilon-\mu-\sqrt{\varepsilon}F_1}{\sqrt{\varepsilon}\sqrt{\mu}}\right\|_2+\sqrt{\varepsilon}^3\left\|w_l \frac{F_\varepsilon-\mu-\sqrt{\varepsilon}F_1}{\sqrt{\varepsilon}\sqrt{\mu}}\right\|_\infty\right)\leq C{\varepsilon}^\frac{1}{4}\to 0,
  \end{equation}
  where $F_1$ is defined in \eqref{1.4.1.5} with \eqref{1.4.1.6} and $(\rho_1,u_1,\theta_1)$ satisfies the acoustic system \eqref{1.4.1.10} with the boundary condition \eqref{1.4.3.2}.
  Therefore, we have justified the hydrodynamic limit from the scaled Boltzmann equation with Maxwell reflection boundary condition with the accommodation coefficient $\alpha\in O(1)$ to the acoustic system with slip boundary condition for the half-space problem. We remark that the convergence rate $\varepsilon^{\frac{1}{4}}$ is optimal. This is due to the appearence of the viscous layer at $O(\varepsilon^{\frac{1}{2}})$ with thickness $O(\varepsilon^{\frac{1}{2}})$, whose $L^2$-norm is naturally $O(\varepsilon^{\frac{1}{4}})$.
\end{remark}

\begin{remark}
  In the current paper, our primary focus is to derive boundary conditions for fluid and Prandtl-type equations by solving the Knudsen layer equation with the accommodation coefficient $\alpha \in(0,1]$, which is quite different with the specular and almost specular reflection cases. So for technical simplicity, we only consider the Boltzmann equation with hard sphere collisions here. In fact, employing the techniques developed in Guo-Jang-Jiang \cite{guo2010acoustic} and the recent work Jiang-Luo \cite{jiang2024knudsenboundarylayerequations}, it is quite standard to address all ranges of collision kernels.  The proofs for general collisions are similar but intricate, so we omit the details here.
\end{remark}

\section{Formal analysis}\label{Se1.4}
In this section, we aim to derive the formal analysis using the Hilbert expansion method. We illustrate how and why boundary layers appear based on the theory developed in Sone's book \cite{sone2002kinetic,sone2007molecular}. Notice that the thickness of the Knudsen layer is $O(\varepsilon)$ and the thickness of the viscous layer is $O(\sqrt{\varepsilon})$.

\subsection{Hilbert expansion}
We seek a solution $F_\varepsilon$ that changes moderately in both $x$ and $t$,
\begin{equation}\label{1.4.1.1}
  F_\varepsilon(t,x,v)=F_0(t,x,v)+\sqrt{\varepsilon}F_1(t,x,v)+\sqrt{\varepsilon}^2F_2(t,x,v)+\cdots,
\end{equation}
where $F_0=\mu$ is the global Maxwellian defined in \eqref{1.3.2}. Substituting the expansion \eqref{1.4.1.1} into the scaled Boltzmann equation \eqref{1.1.1} and arranging the terms by the order of $\varepsilon$, a series of equation for $F_k$ is obtained as follows:
\begin{equation}\label{1.4.1.2}
   \sqrt{\varepsilon}^{-1}: \quad0=Q(\mu,F_1)+Q(F_1,\mu),
  \end{equation}
   \begin{equation}\label{1.4.1.3}
    \sqrt{\varepsilon}^{0}:\quad0=Q(\mu,F_2)+Q(F_2,\mu)+Q(F_1,F_1),
  \end{equation}
\begin{equation}\label{1.4.1.4}
    \sqrt{\varepsilon}^{k}:\quad (\partial_t+v\cdot\nabla_x)F_k=\sum_{i+j=k+2}Q(F_i,F_j), \quad (k\geq 1).
\end{equation}
We define
\begin{equation}\label{1.4.1.5}
  f_k:=\frac{F_k}{\sqrt{\mu}}.
\end{equation}
The solutions of \eqref{1.4.1.2} and \eqref{1.4.1.3} are given as

\begin{equation*}
 \mathcal{L}f_1=0,
 \end{equation*}
 and
 \begin{equation*}
\mathcal{L}f_2=\Gamma(f_1,f_1).
 \end{equation*}
That is,
\begin{equation}\label{1.4.1.6}
  f_1=\left[\rho_1+u_1 \cdot v+\theta_1\left(\frac{|v|^2}{2}-\frac{3}{2}\right)\right]\sqrt{\mu},
\end{equation}
and
\begin{equation}\label{1.4.1.7}
  f_2=\left[\rho_2+u_2 \cdot v+\theta_2\left(\frac{|v|^2}{2}-\frac{3}{2}\right)\right]\sqrt{\mu}+\frac{1}{2}\left[A:u_1\otimes u_1+2\theta_1u_1 \cdot B(v)+\theta_1^2C(v)\right],
\end{equation}
where $A\in \mathbb{R}^{3\times 3},B\in \mathbb{R}^3$ and $C\in\mathbb{R}$ are the Burnett functions defined in Appendix \ref{App.A} and  we used the fact that(see \cite[p.649]{guo2006boltzmann})
\begin{equation}\label{1.4.1.8}
  \mathcal{L}^{-1}\{\Gamma(f,g)+\Gamma(g,f)\}=(\mathbf{I-P})\left\{\frac{\mathbf{P}f\cdot\mathbf{P}g}{\sqrt{\mu}}\right\}.
\end{equation}
Equation \eqref{1.4.1.4} for $f_k$ $(k \geq 1)$ can be solved provided that the solvability condition:

\begin{equation}\label{1.4.1.9}
  \int (\partial_t+v\cdot \nabla_x)f_k\chi_i dv=0, \quad(i=0,1,2,3,4).
  \end{equation}
  where $\chi_i\,(i=0,1,2,3,4)$ are defined in \eqref{1.3.4}.

Substituting \eqref{1.4.1.6} into the solvability condition \eqref{1.4.1.9} for $k=1,$ we obtain the following acoustic system:

\begin{equation}\label{1.4.1.10}
  \begin{cases}
    \partial_t \rho_1+\nabla_x\cdot u_1=0,\\
    \partial_t u_1+\nabla_x(\rho_1+\theta_1)=0,\\
    \partial_t \theta_1+\frac{2}{3}\nabla_x\cdot u_1=0.
  \end{cases}
\end{equation}
Now we consider $k\geq 2$. We decompose $f_k$ into the macroscopic and microscopic parts as follows:
\begin{equation}\label{1.4.1.11}
  f_k=\left[\rho_k+u_kv+\theta_k\left(\frac{|v|^2-3}{2}\right)\right]\sqrt{\mu}+\mathbf{(I-P)}f_k.
\end{equation}
Combining \eqref{1.4.1.9} and \eqref{1.4.1.11}, we get the system for the fluid variables $(\rho_k,u_k,\theta_k)$:
\begin{equation}\label{1.4.1.12}
  \begin{cases}
    \partial_t \rho_k+\nabla_x\cdot u_k=0,\\
    \partial_t u_k+\nabla_x(\rho_k+\theta_k)=- \nabla_x\cdot\int_{\mathbb{R}^3}A \mathbf{(I-P)}f_kdv,\\
    \partial_t \theta_k+\frac{2}{3}\nabla_x\cdot u_k=-\frac{2}{3} \nabla_x\cdot\int_{\mathbb{R}^3}B \mathbf{(I-P)}f_kdv.
  \end{cases}
\end{equation}
Furthermore, the initial conditions of \eqref{1.4.1.12} are imposed on
\begin{equation}\label{1.4.1.13}
  (\rho_k,u_k,\theta_k)(0,x)=(\rho^{in}_k,u^{in}_k,\theta^{in}_k)\in \mathbb{R}\times \mathbb{R}^3\times \mathbb{R}, \quad k=1,2,3\cdots.
\end{equation}
From \eqref{1.4.1.4}, the microscopic part $(\mathbf{I-P})f_k$ on the RHS of \eqref{1.4.1.11} can be determined by $f_i$ for $1\leq i\leq k-1$ as follows:
\begin{equation}\label{1.4.1.14}
  \mathbf{(I-P)}f_k=\mathcal{L}^{-1}\left(\sum_{\stackrel{i+j=k}{i,j\geq1}} \Gamma(f_i,f_{j})-(\partial_t+v\cdot \nabla_x)f_{k-2}\right).
\end{equation}

\subsection{Viscous boundary Layer}
Obviously, $F_1$ in \eqref{1.4.1.6} with \eqref{1.4.1.5} matches the Maxwell reflection boundary condition \eqref{1.1.3} if and only if
\begin{equation}\label{s1}
  u^0_{1}=0,\quad\theta^0_1=0.
\end{equation}
However, the constraints \eqref{s1} on the boundary  are too many as the boundary condition for the acoustic system \eqref{1.4.1.10}. Therefore, this approach does not work and the boundary-layer correction is naturally needed. We express the solution $F_\varepsilon+F^b_\varepsilon$ in the layer with the thickness of $O(\sqrt{\varepsilon})$ (see \cite{sone2002kinetic,sone2007molecular} for example). We first introduce the scaled normal coordinate:
\begin{equation}\label{1.4.2.1}
  \zeta=\frac{x_3}{\sqrt{\varepsilon}}.
\end{equation}
We expand $F_\varepsilon^b$ in a power series of $\sqrt{\varepsilon}$ as
\begin{equation}\label{1.4.2.2}
  F^b_\varepsilon(t,\bar{x},\zeta,v)=\sqrt{\varepsilon}F^b_1(t,\bar{x},\zeta,v)+\sqrt{\varepsilon}^2F^b_2(t,\bar{x},\zeta,v)+\cdots.
\end{equation}
Plugging $F_\varepsilon+F^b_\varepsilon$ into the Boltzmann equation \eqref{1.1.1} gives

  \begin{equation}\label{1.4.2.3}
    \sqrt{\varepsilon}^{-1}:  0=Q\left(\mu, F_{1}^{b}\right)+Q\left(F_{1}^{b}, \mu\right), \\
  \end{equation}

  \begin{equation}\label{1.4.2.4}
\begin{aligned}
      \sqrt{\varepsilon}^{0}:  v_{3} \cdot \partial_{\zeta} F_{1}^{b}=&\left[Q\left(\mu, F_{2}^{b}\right)+Q\left(F_{2}^{b}, \mu\right)\right]+\left[Q\left(F_{1}^{0}, F_{1}^{b}\right)+Q\left(F_{1}^{b}, F_{1}^{0}\right)\right] \\
       &+Q\left(F_{1}^{b}, F_{1}^{b}\right),
\end{aligned}
  \end{equation}

\begin{equation}\label{1.4.2.5}
  \begin{aligned}
  \sqrt{\varepsilon}: \quad &\partial_{t} F_{1}^{b}+\bar{v} \cdot \nabla_{\bar{x}} F_{1}^{b}+v_{3} \cdot \partial_{\zeta} F_{2}^{b}\\
  =&\left[Q\left(\mu, F_{3}^{b}\right)+Q\left(F_{3}^{b}, \mu\right)\right]+\left[Q\left(F_{1}^{0}, F_{2}^{b}\right)+Q\left(F_{2}^{b}, F_{1}^{0}\right)\right] \\
  &+\left[Q\left(F_{2}^{0}, F_{1}^{b}\right)+Q\left(F_{1}^{b}, F_{2}^{0}\right)\right] +\left[Q\left(F_{1}^{b}, F_{2}^{b}\right)+Q\left(F_{2}^{b}, F_{1}^{b}\right)\right]\\
  &+\frac{\zeta}{1 !}\left[Q\left(F_{1}^{(1)}, F_{1}^{b}\right)+Q\left(F_{1}^{b}, F_{1}^{(1)}\right)\right] ,
\end{aligned}
\end{equation}

  \begin{equation}\label{1.4.2.6}
    \begin{aligned}
  \sqrt{\varepsilon}^{k}: \quad \partial_{t} F_{k}^{b}&+\bar{v} \cdot \nabla_{\bar{x}} F_{k}^{b}+v_{3} \cdot \partial_{\zeta} F_{k+1}^{b}
  =\left[Q\left(\mu, F_{k+2}^{b}\right)+Q\left(F_{k+2}^{b}, \mu\right)\right] \\
  &+\sum_{\stackrel{i+j=k+2,}{i, j \geq 1}} Q\left(F_{i}^{b}, F_{j}^{b}\right)+\sum_{\stackrel{i+j=k+2, }{i, j \geq 1}}\left[Q\left(F_{i}^{0}, F_{j}^{b}\right)+Q\left(F_{j}^{b}, F_{i}^{0}\right)\right] \\
 & +\sum_{\stackrel{i+j+l=k+2,}{ 1 \leq l \leq N, i, j \geq 1}} \frac{\zeta^{l}}{l!}\left[Q\left(F_{i}^{(l)}, F_{j}^{b}\right)+Q\left(F_{j}^{b}, F_{i}^{(l)}\right)\right] \quad \text { for } k \geq 1, \\
\end{aligned}
\end{equation}
where we have used the Taylor expansion at $x_3=0$:
\begin{equation*}
  F_i(t,\bar{x},\zeta,v)=F_i^0+\sum_{1\leq l\leq N}\frac{\zeta^l}{l!}F^{(l)}_i+\frac{\zeta^N}{(N+1)!}\tilde{F}_i^{(N+1)}, \quad i\geq1.
\end{equation*}
We define \begin{equation}\label{1.4.2.7}
  f^b_k:=\frac{F^b_k}{\sqrt{\mu}},
\end{equation}
which can be decomposed as
\begin{equation}
  f^b_k=\left[\rho^b_k+u^b_kv+\theta^b_k\left(\frac{|v|^2-3}{2}\right)\right]\sqrt{\mu}+(\mathbf{I-P})f_k^b.
\end{equation}
Throughout this paper, we always assume the following condition
\begin{equation}\label{1.4.2.8}
  f^b_k\to 0,\quad \text{ as } \zeta \to \infty.
\end{equation}
Direct calculations show that the solution of \eqref{1.4.2.3} and \eqref{1.4.2.4} are given as
\begin{equation}\label{1.4.2.9}
  f^b_1=\left[\rho^b_1+u^b_1v+\theta^b_1\left(\frac{|v|^2}{2}-\frac{3}{2}\right)\right]\sqrt{\mu},
\end{equation}
and
\begin{multline}\label{1.4.2.30}
    f^b_2=\left[\rho^b_2+u^b_2v+\theta^b_2\left(\frac{|v|^2-3}{2}\right)\right]\sqrt{\mu}-\sum_{i=1}^{2}\hat{A}_{i3}\partial_\zeta u^b_{1,i}-\hat{B}_3\partial_\zeta\theta^b_1\\
    +\mathcal{L}^{-1}\{\Gamma(f^0_1,f^b_1)+\Gamma(f^b_1,f^0_1)+\Gamma(f^b_1,f^b_1)\}.
\end{multline}

We denote
\begin{equation}
  p^b_k=\rho^b_k+\theta^b_k.
\end{equation}
Multiplying \eqref{1.4.2.4} by $1,v_3 $ and integrating over $v\in \mathbb{R}^3$, one has
\begin{equation*}
  \partial_\zeta u^b_{1,3}=0,\quad \partial_\zeta p^b_1=0,
\end{equation*}
which, together with the far field condition \eqref{1.4.2.8}, one obtains
\begin{equation}\label{1.4.2.10}
  u^b_{1,3}(t,\bar{x},\zeta)\equiv0,\quad  p^b_1(t,\bar{x}, \zeta)\equiv 0, \quad \forall (t,\bar{x},\zeta)\in \mathbb{R}_+\times \mathbb{R}^2\times \mathbb{R}_+.
\end{equation}
$f^b_k$ can be constructed inductively as follows:
\begin{lemma}\label{Lm1.1}
  Let $f^b_k$ be the solution of \eqref{1.4.2.6}, then $(u^b_{k,1},u^b_{k,2},\theta^b_k)$ $(k\geq 1)$ satisfy the following linear Prandtl-type equations:
  \begin{equation}\label{1.4.2.11}
    \begin{cases}
      \partial_tu^b_{k,i}-\kappa_1\partial^2_\zeta u^b_{k,i}+\partial_\zeta\left[(u^b_{1,3}+u^0_{1,3})u^b_{k,i}\right]=\mathrm {f}^b_{k-1,i} \quad (i=1,2),\\
        \partial_t\theta^b_k-\frac{2}{5}\kappa_2\partial^2_\zeta\theta^b_k+\partial_\zeta\left[(u^b_{1,3}+u^0_{1,3})\theta^b_k\right]:=\mathrm{g}^b_{k-1},
    \end{cases}
  \end{equation}
where $\kappa_1,\kappa_2$ are positive constants defined in Appendix \ref{App.A}, and $\mathrm{f}^b_{k-1,i}$ and $g^b_{k-1}$ are defined in \eqref{1.4.2.17} and \eqref{1.4.2.21}, respectively.  Once we have solved $(u^b_{k,1},u^b_{k,2},\theta^b_k)(k\geq 1)$, then $(\mathbf{I-P})f^b_{k+1},u^b_{k+1,3},p^b_{k+1}$ can be determined by \eqref{1.4.2.13}, \eqref{1.4.2.23} and \eqref{1.4.2.19}, respectively.
\end{lemma}

\begin{proof}
  Multiplying \eqref{1.4.2.6} by $1,v,\frac{|v|^2-3}{2}$ and integrating over $\mathbb{R}^3$ to obtain
  \begin{equation}\label{1.4.2.12}
    \begin{cases}
      \partial_t\rho^b_k+\sum_{i=1}^2\partial_iu_{k,i}^b+\partial_\zeta u^b_{k+1,3}=0,\\
      \partial_tu^b_{k,i}+\partial_i(\rho^b_k+\theta^b_k)+\int_{\mathbb{R}^3}A_{i3}(\mathbf{I-P})\partial_\zeta f^b_{k+1}dv+\sum_{j=1}^{2}\int_{\mathbb{R}^3}A_{ij}(\mathbf{I-P})\partial_jf^b_{k}dv=0   \quad (i=1,2),\\
      \partial_tu^b_{k,3}+\partial_\zeta(\rho^b_{k+1}+\theta^b_{k+1})+\int_{\mathbb{R}^3}A_{33}(\mathbf{I-P})\partial_\zeta f^b_{k+1}dv+\sum_{j=1}^{2}\int_{\mathbb{R}^3}A_{3j}(\mathbf{I-P})\partial_jf^b_{k}dv=0,   \\
      \frac{3}{2}\partial_t\theta^b_k-\partial_t\rho^b_k+\int_{\mathbb{R}^3}B_3\mathbf{(I-P)}\partial_\zeta f^b_{k+1}dv+\sum_{j=1}^{2}\int_{\mathbb{R}^3}B_j\mathbf{(I-P)}\partial_jf^b_kdv=0.
    \end{cases}
  \end{equation}
  It is known from \eqref{1.4.2.6} that
  \begin{equation}\label{1.4.2.13}
    \begin{aligned}
      \mathcal{L}f^b_{k+1}=&-\mathbf{(I-P)}(\partial_t+\bar{v}\cdot \nabla_{\bar{x}})f^b_{k-1}-(\mathbf{I-P})[v_3\partial_\zeta f^b_k]\\
      &+\sum_{\stackrel{i+j=k+1,}{i, j \geq 1}} \Gamma\left(f_{i}^{b}, f_{j}^{b}\right)+\sum_{\stackrel{i+j=k+1, }{i, j \geq 1}}\left[\Gamma\left(f_{i}^{0}, f_{j}^{b}\right)+\Gamma\left(f_{j}^{b}, f_{i}^{0}\right)\right] \\
      & +\sum_{\stackrel{i+j+l=k+1,}{ 1 \leq l \leq N, i, j \geq 1}} \frac{\zeta^{l}}{l!}\left[\Gamma\left(f_{i}^{(l)}, f_{j}^{b}\right)+\Gamma\left(f_{j}^{b}, f_{i}^{(l)}\right)\right] \\
      =&-\left(\sum_{j=1}^{3}A_{3j}\partial_\zeta u^b_{k,j} +B_3\partial_\zeta \theta^b_k\right)+\Gamma(f_1^0,\mathbf{P}f^b_k)+\Gamma(\mathbf{P}f^b_k,f^0_1)+\Gamma(f_1^b,\mathbf{P}f^b_k)\\
      &+\Gamma(\mathbf{P}f^b_k,f^b_1)+\mathcal{L}J_{k-1},
    \end{aligned}
  \end{equation}
  where
  \begin{equation} \label{1.4.2.14}
    \begin{aligned}
     J_{k-1}=&\mathcal{L}^{-1}\left\{{-(\mathbf{I-P})[v_3(\mathbf{I-P})\partial_\zeta f^b_k]-(\mathbf{I-P}) (\partial_t+\bar{v}\cdot\nabla_{\bar{x}})f^b_{k-1}} \right.\\
     &\phantom{=\;\;}
     \left.{+\sum_{\stackrel{i+j=k+1, }{i\geq2, j \geq 1}}\left[\Gamma\left(f_{i}^{0}, f_{j}^{b}\right)+\Gamma\left(f_{j}^{b}, f_{i}^{0}\right)\right]
   +\sum_{\stackrel{i+j=k+1,}{i, j \geq 2}} \Gamma\left(f_{i}^{b}, f_{j}^{b}\right)} \right.\\
   &\phantom{=\;\;}
   \left.\sum_{\stackrel{i+j+l=k+1, }{i,j \geq1, 1\leq l \leq N}}\frac{\zeta^l}{l!}\left[\Gamma(f^{(l)}_i,f^b_j)+\Gamma(f^b_j,f^{(l)}_i)\right]\right.\\
   &\phantom{=\;\;}
   \left.{
   +\Gamma(f^0_1,(\mathbf{I-P})f^b_k)+\Gamma((\mathbf{I-P})f^b_k,f^0_1) +\Gamma(f^b_1,(\mathbf{I-P})f^b_k)+\Gamma((\mathbf{I-P})f^b_k,f^b_1)}
   \right\}.
  \end{aligned}
  \end{equation}
  Using \eqref{1.4.1.6}, \eqref{1.4.1.8} and \eqref{1.4.2.13}, one obtains
  \begin{equation}\label{1.4.2.15}
    \begin{aligned}
      \mathbf{(I-P)} f^b_{k+1}=&-\left(\sum_{j=1}^{3}\hat{A}_{3j}\partial_\zeta u^b_{k,j} +\hat{B}_3\partial_\zeta \theta^b_k\right)\\
       &+\left[A(v):u^b_1\otimes u^b_k+B(v)(u^b_1\theta^b_k+u^b_k\theta^b_1)+C(v)\theta^b_1\theta^b_k\right]\\
       &+\left[A(v):u^0_1\otimes u^b_k+B(v)(u^0_1\theta^b_k+u^b_k\theta^0_1)+C(v)\theta^0_1\theta^b_k\right]+J_{k-1}.
    \end{aligned}
  \end{equation}
  Substituting \eqref{1.4.2.15} into $\eqref{1.4.2.12}_2$, one obtains that

  \begin{equation}\label{1.4.2.16}
    \partial_tu^b_{k,i}-\kappa_1\partial^2_\zeta u^b_{k,i}+\partial_\zeta\left[(u^b_{1,3}+u^0_{1,3})u^b_{k,i}\right]=\mathrm {f}^b_{k-1,i} \quad (i=1,2),
  \end{equation}
  where we have used(see \cite[lemma 4.4]{bardos1993fluid} for instance)
  \begin{equation*}
    \quad \left<\hat{A}_{ij},A_{kl}\right>=\kappa_1(\delta_{ij}\delta_{jl}+\delta_{il}\delta_{jk}-\frac{2}{3}\delta_{ij}\delta_{kl}).
  \end{equation*}
  The source term $\mathrm{f}^b_{k-1,i}$ is given as
  \begin{equation}\label{1.4.2.17}
    \mathrm {f}^b_{k-1,i}=-\partial_\zeta\left[(u^b_{1,i}+u^0_{1,i})u^b_{k,3}\right]-\left<A_{3i},\partial_\zeta J_{k-1}\right>+W^b_{k-1,i}-\partial_ip^b_k,
  \end{equation}
  and
  \begin{equation}\label{1.4.2.18}
    W^b_{k-1,i}=-\sum_{j=1}^{2}\left<A_{ij}\partial_j(\mathbf{I-P})f^b_k\right> \quad (i=1,2,3).
  \end{equation}
  Using \eqref{1.4.2.15} and $\eqref{1.4.2.12}_3$, one has
\begin{multline}\label{1.4.2.19}
      \partial_\zeta p^b_{k+1}=-\partial_t u^b_{k,3}+\frac{4}{3}\kappa_1\partial^2_\zeta u^b_{k,3}+\frac{2}{3}\sum_{i=1}^{2}\partial_\zeta\left[(u^b_{1,i}+u^0_{1,i})u^b_{k,i}\right]\\
      -\frac{4}{3}\partial_\zeta\left[(u^b_{1,3}+u^0_{1,3})u^b_{k,3}\right]
      -\left<A_{33},\partial_\zeta J_{k-1}\right>+W^b_{k-1,3},
\end{multline}
  where we have used the fact that
  \begin{equation*}
    \left<A_{ii},A_{jj}\right>=-\frac{2}{3},\quad \left<A_{ij},A_{ij}\right>=1, \text{ for } i\neq j, \quad \text{and }\left<A_{ii},A_{ii}\right>=\frac{4}{3}.
  \end{equation*}
  Similarly, from \eqref{1.4.2.15} and $\eqref{1.4.2.12}_4$, the equation for $\theta^b_k$ can be reduced to be
  \begin{equation}\label{1.4.2.20}
    \partial_t\theta^b_k-\frac{2}{5}\kappa_2\partial^2_\zeta\theta^b_k+\partial_\zeta\left[(u^b_{1,3}+u^0_{1,3})\theta^b_k\right]:=\mathrm{g}^b_{k-1},
  \end{equation}
  with
  \begin{equation}\label{1.4.2.21}
    \mathrm{g}^b_{k-1}=-\partial_\zeta\left[(\theta^b_1+\theta^0_1)u^b_{k,3}\right]-\frac{2}{5}\left<B_3,\partial_\zeta J_{k-1}\right>+\frac{2}{5}H_{k-1}+\frac{2}{5}\partial_tp^b_k,
  \end{equation}
  and
  \begin{equation}\label{1.4.2.22}
    H^b_{k-1}=-\sum_{j=1}^{2}\left<B_i,(\mathbf{I-P})\partial_if^b_k\right>.
  \end{equation}
  Moreover, $\eqref{1.4.2.12}_1$ gives
  \begin{equation}\label{1.4.2.23}
    \partial_\zeta u^b_{k+1,3}=-\partial_t\rho^b_k-\sum_{i=1}^2\partial_iu_{k,i}^b.
  \end{equation}
\end{proof}

We point out that $W^b_{k-1},H^b_{k-1}$ and $J^b_{k-1}$ depend on $f^b_j,1\leq j\leq k-1$. Moreover, $J^b_0=W^b_0=H^b_0=\mathrm{f}^b_0=\mathrm{g}^b_0=0.$ Finally, we impose the initial conditions of \eqref{1.4.2.11} on
\begin{equation}\label{1.4.2.31}
  (\bar{u}^b_k,\theta^b_k)(0,\bar{x},\zeta)=(\bar{u}^{b,in}_k,\theta^{b,in}_k)(\bar{x},\zeta)\in \mathbb{R}^2\times \mathbb{R},\quad k=1,2,3,\cdots,
\end{equation}
with
\begin{equation*}
  \lim_{\zeta\to \infty}(\bar{u}^{b,in}_k,\theta^{b,in}_k)(\bar{x},\zeta)=0.
\end{equation*}

\subsection{Knudsen boundary layer and boundary conditions}\label{se1.4.3}

Since the sum of the Hilbert and viscous boundary-layer solutions is a infinitesimal Maxwellian at the first order of $\sqrt{\varepsilon}$, it matches the Maxwell reflection boundary condition at this order provided that
\begin{equation}\label{1.4.3.1}
\begin{array}{l}
    u^0_{1,3}+u^{b,0}_{1,3}=0, \\
    u^{b,0}_{1,i}+u^0_{1,i}=0,\quad (i=1,2),\\
    \theta^{b,0}_1+\theta^0_1=0,\\
\end{array}
\end{equation}
are satisfied on the boundary. Here the subscript $0$ indicates the value on the boundary, i.e., $x_3=0$. Therefore, the Knudsen-layer correction is not required at this order. Utilizing \eqref{1.4.2.10} and \eqref{1.4.3.1}, we obtain the boundary condition for \eqref{1.4.1.10} :
\begin{equation}\label{1.4.3.2}
  u^0_{1,3}=0,
\end{equation}
and the boundary condition for \eqref{1.4.2.11} with $k=1$:
\begin{equation}\label{1.4.3.3}
\begin{aligned}
   & u^{b,0}_{1,i}=-u^0_{1,i}, \quad( i=1,2),\\
   & \theta^{b,0}_1=-\theta^0_1.
\end{aligned}
\end{equation}

Remark that, the boundary condition of the acoustic system \eqref{1.4.3.2} coincide with the boundary condition derived in \cite{jiang2010remarks}. Moreover, by the boundary conditions derived above, we have the following remark:
\begin{remark}\label{Rm1.2}
  We remark that by using \eqref{1.4.2.10} and $\eqref{1.4.3.2}$, the Prandtl-type equations \eqref{1.4.2.11} can be reduced to be the following heat equations:
  \begin{equation}\label{1.4.3.4}
    \begin{cases}
      \partial_tu^b_{k,i}-\kappa_1\partial^2_\zeta u^b_{k,i} =\mathrm {f}^b_{k-1,i} \quad (i=1,2),\\
        \partial_t\theta^b_k-\frac{2}{5}\kappa_2\partial^2_\zeta\theta^b_k=\mathrm{g}^b_{k-1},
    \end{cases}
  \end{equation}
\end{remark}

However, the solution $F_\varepsilon+F^b_\varepsilon$ does not match the boundary condition at $O(\sqrt{\varepsilon}^k)$ for $(k\geq2)$ and the Knudsen-layer correction $F_\varepsilon^{bb}$ is required. From \cite{sone2002kinetic,sone2007molecular}, the Knudsen-layer appears only adjacent to the boundary with thickness $O(\varepsilon)$ and vanishes sufficiently fast away from the boundary. Inside the Knudsen layer, the new stretched coordinate is needed:
\begin{equation}\label{1.4.3.5}
  \xi=\frac{x_3}{\varepsilon}.
\end{equation}
Taking into account that the correction is not required at the leading and the first order of $\sqrt{\varepsilon}$, we define the Knudsen boundary expansion as:
\begin{equation}\label{1.4.3.6}
  F^{bb}_\varepsilon(t,\bar{x},\xi,v)=\sqrt{\varepsilon}^2F^{bb}_2(t,\bar{x},\xi,v)+\sqrt{\varepsilon}^3F^{bb}_3(t,\bar{x},\xi,v)+\cdots,
\end{equation}
 Plugging $F_\varepsilon+F_\varepsilon^b+F_\varepsilon^{bb}$ into the Boltzmann equation \eqref{1.1.1}, one obtains that

\begin{equation}\label{1.4.3.7}
  \sqrt{\varepsilon}^0: \quad v_3\partial_\xi F^{bb}_2 =Q( F^{bb}_2,\mu)+Q(\mu,F^{bb}_2),
\end{equation}
\begin{equation}\label{1.4.3.8}
  \sqrt{\varepsilon}^1: \quad v_3\partial_\xi F^{bb}_3= Q(\mu, F^{bb}_3)+Q( F^{bb}_3,\mu)+Q( F^{bb}_2, F^0_1)+Q( F^0_1, F^{bb}_2)\\
 +Q( F^{bb}_2, F^{b,0}_1)+Q( F^{b,0}_1,F^{bb}_2),
\end{equation}

\begin{equation}\label{1.4.3.9}
  \begin{aligned}
\sqrt{\varepsilon}^{k}: \quad \partial_{t} F_{k}^{bb}&+\bar{v} \cdot \nabla_{\bar{x}} F_{k}^{bb}+v_{3} \cdot \partial_{\xi} F_{k+2}^{bb}
=\left[Q\left(\mu, F_{k+2}^{bb}\right)+Q\left(F_{k+2}^{bb}, \mu\right)\right] \\
&+\sum_{\stackrel{i+j=k+2,}{i, j \geq 1}} Q\left(F_{i}^{bb}, F_{j}^{0}+F^{b,0}_j\right)+Q\left(F_{j}^{0}+F^{b,0}_j,F_{i}^{bb}\right)+\sum_{\stackrel{i+j=k+2, }{i, j \geq 1}} Q\left(F_{i}^{bb}, F_{j}^{bb}\right) \\
& +\sum_{\stackrel{i+j+2l=k+2,}{ 1 \leq l \leq N, i\geq 2, j \geq 1}} \frac{\xi^{l}}{l!}\left[Q\left(F_{i}^{bb}, F_{j}^{(l)}\right)+Q\left(F_{j}^{(l)}, F_{i}^{bb}\right)\right]\\
& +\sum_{\stackrel{i+j+l=k+2,}{ 1 \leq l \leq N, i\geq 2, j \geq 1}} \frac{\xi^{l}}{l!}\left[Q\left(F_{i}^{bb}, F_{j}^{b,(l)}\right)+Q\left(F_{j}^{b,(l)}, F_{i}^{bb}\right)\right] \quad \text { for } k \geq 1, \\
\end{aligned}
\end{equation}
where the Taylor expansion of $F^{b}_i$ at $\zeta=0$ is utilized:
\begin{equation*}
  F^b_i(t,\bar{x},\xi,v)=F_i^{b,0}+\sum_{1\leq l\leq N}\frac{\xi^l}{l!}F^{b,(l)}_i+\frac{\xi^N}{(N+1)!}\tilde{F}_i^{b,(N+1)}, \quad i\geq1.
\end{equation*}
Similar in the viscous layer, we define $f^{bb}_k=\frac{F^{bb}_k}{\sqrt{\mu}}.$ We also assume the following far field condition
\begin{equation}\label{1.4.3.10}
  f^{bb}_k\to 0, \quad \text{as } \xi\to \infty.
\end{equation}
Collecting \eqref{1.1.3},\eqref{1.4.3.9} and \eqref{1.4.3.10}, we then obtain the following boundary layer problem for $k\geq 2$
\begin{equation}\label{1.4.3.11}
\begin{cases}
    v_3\partial_\xi f^{bb}_k+\mathcal{L}f^{bb}_k=S^{bb}_k\quad (\xi>0),\\
    f^{bb}_k=\mathcal{K}f^{bb}_k-\tilde{\mathcal{K}}(f_k+f^b_k)   \quad (\xi=0,v_3>0),\\
    f^{bb}_k \to 0 \quad \text{as } \xi \to \infty .
\end{cases}
\end{equation}
where the boundary operator
\begin{equation}\label{1.4.3.20}
  \mathcal{K}f=(1-\alpha)f(R_xv)+\alpha\sqrt{2\pi}\sqrt{\mu(v)}\int_{v'\cdot n>0}(v'\cdot n)f\sqrt{\mu(v')}dv',
\end{equation}
and $\tilde{\mathcal{K}}f=f-\mathcal{K}f$ ($\tilde{\mathcal{K}}_1$ that appears later is defined in the same way). The source term $S^{bb}_k$ can be divided into $S^{bb}_k=S^{bb}_{k,1}+S^{bb}_{k,2}$ with
\begin{equation}\label{Sbb-1}
  S^{bb}_{k,1}=-\mathbf{P}\{(\partial_t+\bar{v}\cdot\nabla_{\bar{x}})f^{bb}_{k-2}\}\in \mathcal{N},
\end{equation}
\begin{equation}\label{Sbb-2}
\begin{aligned}
    S^{bb}_{k,2}=&-\mathbf{(I-P)}\{(\partial_t+\bar{v}\cdot\nabla_{\bar{x}})f^{bb}_{k-2}\}+\sum_{\stackrel{i+j=k+2,}{i, j \geq 1}} \Gamma\left(f_{i}^{bb}, f_{j}^{0}+f^{b,0}_j\right)\\
    &+\Gamma\left(f_{j}^{0}+f^{b,0}_j,f_{i}^{bb}\right)+\sum_{\stackrel{i+j=k+2, }{i, j \geq 1}} \Gamma\left(f_{i}^{bb}, f_{j}^{bb}\right) \\
    & +\sum_{\stackrel{i+j+2l=k+2,}{ 1 \leq l \leq N, i\geq 2, j \geq 1}} \frac{\xi^{l}}{l!}\left[\Gamma\left(f_{i}^{bb}, f_{j}^{(l)}\right)+\Gamma\left(f_{j}^{(l)}, f_{i}^{bb}\right)\right]\\
    & +\sum_{\stackrel{i+j+l=k+2,}{ 1 \leq l \leq N, i\geq 2, j \geq 1}} \frac{\xi^{l}}{l!}\left[\Gamma\left(f_{i}^{bb}, f_{j}^{b,(l)}\right)+\Gamma\left(f_{j}^{b,(l)}, f_{i}^{bb}\right)\right] \in \mathcal{N}^\perp.\\
\end{aligned}
\end{equation}
Here the notation $f^{bb}_{-1}=f^{bb}_0=f^{bb}_1=0$ are employed. As shown below, the solvability conditions for the Knudsen layer equation \eqref{1.4.3.11} provide the boundary conditions for the macroscopic part of $f_k$ and $f^b_k$, allowing us to solve for $f_k$ and $f^b_k$. Consequently, when solving $f^{bb}_k$, the quantities $S^{bb}_{k,1}$ and $S^{bb}_{k,2}$ are known, as $f_i,f^b_i$ and $f^{bb}_j$ $( i\leq k, j\leq k-1)$ have already been determined.

It should be noted that the boundary condition $u_{k,3}+u^b_{k,3}$ on the boundary $(\xi=0)$ is obtained by multiplying $\eqref{1.4.3.11}_1$ by $\sqrt{\mu}$ and integrating from $\xi=0$ to $+\infty$ under the far-field condition \eqref{1.4.3.10} as
\begin{equation}\label{1.4.3.14}
  u^0_{k,3}+u^{b,0}_{k,3}=-\int_{0}^{+\infty}\int_{\mathbb{R}^3}S^{bb}_k\sqrt{\mu}dvd\xi.
\end{equation}

To solve the Knudsen layer problems and obtain the boundary conditions of the fluid equations, we firstly consider $k=2$. Plugging \eqref{1.4.1.7} and \eqref{1.4.2.30} into \eqref{1.4.3.11}, and using \eqref{1.4.3.7}, the boundary value problem \eqref{1.4.3.11} is reduced to
\begin{equation}\label{1.4.3.15}
  \begin{cases}
      v_3\partial_\xi f^{bb}_2+\mathcal{L}f^{bb}_2=0 \quad (\xi>0),\\
      \begin{aligned}f^{bb}_2=&\mathcal{K}f^{bb}_2-\tilde{\mathcal{K}}_1[\sqrt{\mu}]\sum_{i=1}^{2}(u^0_{2,i}+u^{b,0}_{2,i})v_i-\tilde{\mathcal{K}}[v_3\sqrt{\mu}](u^0_{2,3}+u^{b,0}_{2,3})-\tilde{\mathcal{K}}[\frac{|v|^2-3}{2}\sqrt{\mu}](\theta^0_2+\theta^{b,0}_2)\\
        &+\tilde{\mathcal{K}}_1[\alpha(|v|)v_3\sqrt{\mu}]\sum_{i=1}^{2}\partial_\zeta u^{b,0}_{1,i}v_i+\tilde{\mathcal{K}}[\hat{B}_3]\partial_\zeta\theta^{b,0}_1 \quad (\xi=0,v_3>0),
      \end{aligned}\\
      f^{bb}_2 \to 0 \quad \text{as } \xi \to \infty.
  \end{cases}
  \end{equation}
Similar as the analysis in \cite{aoki2017slip,hattori2019sound,takata2012asymptotic}, by the isotropic property of $\mathcal{L}$ and $\mathcal{K}$ (see Appendix \ref{App.B}), $f^{bb}_2$ is expressed in terms of two fundamental solutions as
\begin{equation}\label{1.4.3.16}
  f^{bb}_2=\sum_{i=1}^{2}\partial_\zeta u^{b,0}_{1,i}v_i\phi^{(1)}_1(\xi,|v|,v_3)+\partial_\zeta \theta^{b,0}_1 \phi^{(0)}_1(\xi,|v|,v_3).
\end{equation}
Here $\phi^{(1)}_1,\phi^{(0)}_1$ are, respectively, the solutions of the following boundary-value problems:
\begin{equation}\label{1.4.3.17}
  \begin{cases}
    v_3\partial_\xi \phi^{(1)}_1+\mathcal{L}_1\phi^{(1)}_1=0,\\
    \phi^{(1)}_1=\mathcal{K}_1\phi^{(1)}_1-\tilde{\mathcal{K}}_1[\sqrt{\mu}]b_1+\tilde{\mathcal{K}}_1[\alpha(|v|)v_3\sqrt{\mu}]  \quad (\xi=0,v_3>0),\\
    \phi^{(1)}_1\to 0 \quad \text{as }\xi\to \infty,
  \end{cases}
\end{equation}
and
\begin{equation}\label{1.4.3.18}
  \begin{cases}
    v_3\partial_\xi \phi^{(0)}_1+\mathcal{L} \phi^{(0)}_1=0,\\
    \phi^{(0)}_1=\mathcal{K} \phi^{(0)}_1-\tilde{\mathcal{K}}[\frac{|v|^2-3}{2}\sqrt{\mu}]c_1+\tilde{\mathcal{K}}[\hat{B}_3]  \quad (\xi=0,v_3>0),\\
    \phi^{(0)}_1\to 0\quad \text{as }\xi\to \infty.
  \end{cases}
\end{equation}
where $\mathcal{L}_1$ and $\mathcal{K}_1$ are defined in Appendix \ref{App.B}. It is known from \cite{bardos1986milne,coron1988classification,Jiang2024boundary} that for each of the problems \eqref{1.4.3.17}  and \eqref{1.4.3.18} has a unique solution only when $b_1$ and $c_1$ take special values. The constants $b_1$ and $c_1$ are the so-called slip coefficients and determined together with the solutions $\phi^{(1)}_1$ and $\phi^{(0)}_1$, respectively.

By the comparison of \eqref{1.4.3.15} with \eqref{1.4.3.16}-\eqref{1.4.3.18} and using \eqref{1.4.3.14}, we obtain the following boundary condition for the fluid equations:
\begin{equation}
\begin{cases}
    u^0_{2,i}+u^{b,0}_{2,i}=b_1\partial_\zeta u^{b,0}_{1,i},\\
    u^0_{2,3}+u^{b,0}_{2,3}=0,\\
    \theta^0_2+\theta^{b,0}_2=c_1 \partial_\zeta \theta^{b,0}_1.
\end{cases}
\end{equation}
That is,
\begin{equation}
  u^0_{2,3}=-u^{b,0}_{2,3}
\end{equation}
 for \eqref{1.4.1.12} with $k=2$ and
 \begin{equation}
  \begin{cases}
   u^{b,0}_{2,i}=-u^0_{2,i}+b_1\partial_\zeta u^{b,0}_{1,i},\\
\theta^{b,0}_2=-\theta^0_2+c_1 \partial_\zeta \theta^{b,0}_1.
\end{cases}
\end{equation}
for \eqref{1.4.2.11} with $k=2$.

We now consider the $\tfrac{3}{2}$th order of the Knudsen number. Plugging $f_3+f^b_3$ into \eqref{1.4.3.11} and using \eqref{1.4.3.8}, one obtains the system for $f^{bb}_3$:

  \begin{equation}
    \begin{cases}
      v_3\partial_\xi f^{bb}_3+\mathcal{L}f^{bb}_3=-(\rho^0_1+\rho^{b,0}_1)\partial_\zeta u^{b,0}_{1,i}v_1\mathcal{L}_1\phi^{(1)}_1-(\rho^0_1+\rho^{b,0}_1)\partial_\zeta \theta^{b,0}_{1}\mathcal{L}\phi^{(0)}_1,\\
        \begin{aligned}
          f^{bb}_3=&\mathcal{K}f^{bb}_3-\sum_{i=1}^{2}(u^0_{3,i}+u^{b,0}_{3,i})v_i\tilde{\mathcal{K}}_1[\sqrt{\mu}]-(u^0_{3,3}+u^{b,0}_{3,3})\tilde{\mathcal{K}}[v_3\sqrt{\mu}]-\tilde{\mathcal{K}}[\frac{|v|^2-3}{2}\sqrt{\mu}](\theta^0_3+\theta^{b,0}_3)\\
          &+ \tilde{\mathcal{K}}_1[v_3\alpha(|v|)\sqrt{\mu}]\sum_{i=1}^{2}v_i\left(\partial_3u^0_{1,i}+\partial_\zeta u^{b,0}_{2,i}-(\rho^0_1+\rho^{b,0}_1)\partial_\zeta u^b_{1,i}\right)\\
          &+ \tilde{\mathcal{K}}[\hat{A}_{33}](\partial_3u^0_{1,3}+\partial_\zeta u^{b,0}_{2,3})\\
          &+\tilde{\mathcal{K}}[\hat{B}_3]\left[\partial_3\theta^0_1+\partial_\zeta\theta^{b,0}_2-(\rho^0_1+\rho^{b,0}_1)\partial_\zeta\theta^{b,0}_1\right]\\
          &-\tilde{\mathcal{K}}_1[D_1(|v|)+v_3^2D_2(|v|)]\sum_{i=1}^{2}v_i\partial_\zeta^2u^{b,0}_{1,i}-\tilde{\mathcal{K}}[F_2(|v|,v_3)]\partial_\zeta^2\theta^{b,0}_1\quad (\xi=0,v_3>0),
        \end{aligned}\\
        f^{bb}_3 \to 0 \quad \text{as } \xi \to \infty,
    \end{cases}
    \end{equation}
where we have used the boundary condition \eqref{1.4.3.1} and $D_1(|v|),D_2(|v|)$ and $F_2(|v|,v_3)$ are only functions of $|v|,v_3$ and defined in Appendix \ref{App.A}. After the same arguments as \eqref{1.4.3.15}-\eqref{1.4.3.18}, we can express $f^{bb}_3$ as
\begin{equation}
  f^{bb}_3=f^{bb,(0)}_3+f^{bb,(1)}_3,
\end{equation}
each part of which is expressed as follows:
 
\begin{multline}
    f^{bb,(0)}_3=\phi^{(0)}_1\left[\partial_3\theta^0_1+\partial_\zeta\theta^{b,0}_2-\partial_\zeta\theta^{b,0}_1(\rho^{b,0}_1+\rho^0_1)\right]\\
    +\phi^{(0)}_2\partial^2_\zeta\theta^{b,0}_1+\phi^{(0)}_3(\partial_3u^0_{1,3}+\partial_\zeta u^{b,0}_{2,3})+\phi^{(0)}_4\partial_\zeta \theta^{b,0}_1(\rho^{b,0}_1+\rho^0_1),
\end{multline}
 
and
\begin{equation}
  f^{bb,(1)}_3=\sum_{i=1}^{2}v_i\left(\phi^{(1)}_1\left[\partial_3u_{1,i}^0+\partial_\zeta u^{b,0}_{1,i}-\partial_\zeta u^{b,0}(\rho^{b,0}_1+\rho^0_1)\right]+\phi^{(1)}_2\partial_\zeta^2u^{b,0}_{1,i}+\phi^{(1)}_3\partial_\zeta u^{b,0}_{1,i}(\rho^{b,0}_1+\rho^0_1)\right),
\end{equation}
where the elements $\phi^{(0)}_1$ and $\phi^{(1)}_1$ are the solutions of system \eqref{1.4.3.18} and \eqref{1.4.3.17}, respectively. The other elements $\phi^{(0)}_j(j=2,3,4)$,$\phi^{(1)}_j(j=2,3)$ are the solution of the following problem, respectively:
\begin{equation}
\begin{cases}
      v_3\partial_\zeta\phi^{(0)}_{j}+\mathcal{L}\phi^{(0)}_{j}=I^{(0)}_j,\\
      \phi^{(0)}_j=-c_j\tilde{\mathcal{K}}[\frac{|v|^2-3}{2}\sqrt{\mu}]+\mathcal{K}[\phi^{(0)}_j]+\tilde{\mathcal{K}}[g^{(0)}_j] \quad (\xi=0,v_3>0),\\
      \phi^{(0)}_j\to 0 \quad \text{as } \xi\to \infty,
\end{cases}
\end{equation}
where
\begin{equation}
\begin{aligned}
    I^{(0)}_2=0 , \quad g^{(0)}_2=-F_2(|v|,v_3),\\
    I^{(0)}_3=0, \quad g^{(0)}_3=-\hat{A}_{33},\\
    I^{(0)}_4=-\mathcal{L}\phi^{(0)}_1(|v|,v_3),\quad  g^{(0)}_4=0,
\end{aligned}
\end{equation}
and
\begin{equation}
\begin{cases}
        v_3\partial_\zeta\phi^{(1)}_{j}+\mathcal{L}_1\phi^{(1)}_{j}=I^{(1)}_j,\\
        \phi^{(1)}_j=-b_j\tilde{\mathcal{K}}_1[\sqrt{\mu}]+\mathcal{K}_1[\phi^{(1)}_j]+\tilde{\mathcal{K}}_1[g^{(1)}_j], \quad (\xi=0,v_3>0),\\
        \phi^{(1)}_j\to 0 \quad \text{as } \xi\to \infty,
\end{cases}
  \end{equation}
  where
  \begin{equation}
  \begin{aligned}
      I^{(1)}_2=0 , \quad g^{(1)}_2=-D_1(|v|)-v^2_3D_2(|v|),\\
      I^{(1)}_3=-\mathcal{L}_1\phi^{(1)}_1, \quad g^{(1)}_3=0.
  \end{aligned}
  \end{equation}
Here $\phi^{(i)}_j(i=0,1;j=1,2...)$ on the right-hand side are fundamental solutions that depends only on $\xi,v_3,|v|$. In view of their boundary-value problems, the boundary conditions are obtained
\begin{equation}
  \begin{cases}
    u^{b,0}_{3,i}+u^0_{3,i}=b_1\left[\partial_3u^0_{1,i}+\partial_\zeta u^{b,0}_{2,i}-\partial_\zeta u^{b,0}_{1,i}(\rho^{b,0}_1+\rho^0_1)\right]+b_2\partial^2_\zeta u^{b,0}_{1,i}+b_3(\rho^0_1+\rho^{b,0}_1)\partial_\zeta u^{b,0}_{1,i},\\
u^{b,0}_{3,3}+u^0_{3,3}=0,\\
\begin{aligned}
  \theta_3^{b,0}+\theta^0_3=&c_1\left[\partial_3\theta^0_1+\partial_\zeta\theta^{b,0}_2-(\rho^0_1+\rho^{b,0}_1)\partial_\zeta\theta^{b,0}_1\right]+c_2\partial_\zeta^2\theta^{b,0}_1+c_3(\partial_3u^0_{1,3}+\partial_\zeta u^{b,0}_{2,3})\\
  &+c_4(\rho^0_1+\rho^{b,0}_1)\partial_\zeta\theta^{b,0}_1,
\end{aligned}
  \end{cases}
\end{equation}
for $i=1,2$.

Similarly, $f^{bb}_k(k\geq 4)$ can be constructed by the same ways of constructing $f^{bb}_2$ and $f^{bb}_3$ by using the symmetry properties of operators $\mathcal{L},\Gamma$ and $\mathcal{K}$ (see Appendix \ref{App.B}). The calculations are teidous but trival, so we will omit the details here for simplicity.

\subsection{Truncations of the Hilbert expansion}\label{subse:truncate}
Motivated by Caflish's work \cite{caflisch1980fluid}, in the present paper we aim to rigorously justify the acoustic limit from the scaled Boltzmann equation by the Hilbert expansion approach. Based on the formal analysis above, we truncate the Hilbert expansion as follows:
  \begin{equation}\label{1.5.1}
        F_\varepsilon=\mu+\sum_{k=1}^{6}\sqrt{\varepsilon}^kF_k(t,x,v)+\sum_{k=1}^{6}\sqrt{\varepsilon}^kF^b_k(t,\bar{x},\zeta,v)
         +\sum_{k=2}^{6}\sqrt{\varepsilon}^kF^{bb}_k(t,\bar{x},\xi,v)+\sqrt{\varepsilon}^5\mu^{\frac{1}{2}}f_{R,\varepsilon}\geq 0,
    \end{equation}
where the global Maxwellian $\mu$ is defined in \eqref{1.3.2}. Moreover,  we choose the number $N\in\mathbb{N}_+$ that appeared in the Taylor expansions before by $N=5$ for the Hilbert expansion \eqref{1.5.1}.

Consequently, plugging \eqref{1.5.1} into the scaled Boltzmann equation \eqref{1.1.1} and utilizing the boundary condition \eqref{1.1.3}, we obtain the following  remainder equation

\begin{equation}\label{1.5.2}
  \partial_tf_{R,\varepsilon} +v\cdot\nabla_xf_{R,\varepsilon}+\frac{1}{\varepsilon}\mathcal{L}f_{R,\varepsilon}=S,
\end{equation}
with the Maxwell reflection boundary condition
\begin{equation}\label{1.5.3}
   f_{R,\varepsilon}=\mathcal{K}f_{R,\varepsilon} ,\quad \text{on }\Sigma_-,
\end{equation}
where the Maxwell reflection boundary operator $\mathcal{K}$ is defined in \eqref{1.4.3.20}, and the source term
\begin{equation}\label{1.5.4}
  S:=S_1+S_2+S_3+S_3+S_4+S_5,
\end{equation}
with
\begin{equation}\label{1.5.5}
  S_1:=\sqrt{\varepsilon}^3\Gamma(f_{R,\varepsilon},f_{R,\varepsilon}),
\end{equation}
\begin{equation}\label{1.5.6}
  S_2:=\frac{1}{\varepsilon}\sum_{i=1}^{6}\left[\Gamma(f_{R,\varepsilon},\sqrt{\varepsilon}^i(f_i+f^b_i+f^{bb}_i))+\Gamma(\sqrt{\varepsilon}^i(f_i+f^b_i+f^{bb}_i),f_{R,\varepsilon})\right],
\end{equation}
\begin{equation}\label{1.5.7}
  S_3:=-(\partial_t+v\cdot\nabla_x)(f_5+\sqrt{\varepsilon}f_6)+\sum_{\stackrel{i+j\geq7,}{1\leq i, j \leq 6}}\sqrt{\varepsilon}^{i+j-7}\Gamma(f_i,f_j),
\end{equation}
\begin{multline}\label{1.5.8}
  \begin{aligned}
  S_4:= & -\left(\partial_{t}+\bar{v} \cdot \nabla_{\bar{x}}\right)\left(f_{5}^{b}+\sqrt{\varepsilon} f_{6}^{b}\right)-v_{3} \partial_{\zeta} f_{6}^{b} \\
    & +\sum_{\substack{i+j \geq 7 \\
    1 \leq i, j \leq 6}} \sqrt{\varepsilon}^{i+j-7}\left[\Gamma\left(f_{i}^{0}, f_{j}^{b}\right)+\Gamma\left(f_{j}^{b}, f_{i}^{0}\right)+\Gamma\left(f_{i}^{b}, f_{j}^{b}\right)\right] \\
    & +\sum_{\substack{i+j+l \geq 7\\
    1\leq i,j\leq 6,1\leq l\leq 5}} \sqrt{\varepsilon}{ }^{i+j+l-7} \frac{\zeta^{l}}{l !}\left[\Gamma\left(f_{i}^{(l)}, f_{j}^{b}\right)+\Gamma\left(f_{j}^{b}, f_{i}^{(l)}\right)\right] \\
    & +\frac{\zeta^{6}}{6 !} \sum_{j=1}^{6} \sqrt{\varepsilon}^{j-1}\left[\Gamma\left(\sum_{i=1}^{6} \sqrt{\varepsilon}^{i} \widetilde{f}_{i}^{(6)}, f_{j}^{b}\right)+\Gamma\left(f_{j}^{b}, \sum_{i=1}^{6} \sqrt{\varepsilon}^{i} \widetilde{f}_{i}^{(6)}\right)\right],
    \end{aligned}
\end{multline}
and
\begin{multline}\label{1.5.9}
  \begin{aligned}
  S_5:=&-\left(\partial_{t}+\bar{v} \cdot \nabla_{\bar{x}}\right)\left(f_{5}^{b b}+\sqrt{\varepsilon} f_{6}^{b b}\right)\\
    &+\sum_{\substack{i+j \geq 7 \\
    1 \leq i, j \leq 6}} \sqrt{\varepsilon}^{i+j-7}\left[\Gamma\left(f_{i}^{0}+f_{i}^{b, 0}, f_{j}^{b b}\right)+\Gamma\left(f_{j}^{b b}, f_{i}^{0}+f_{i}^{b, 0}\right)+\Gamma\left(f_{i}^{b b}, f_{j}^{b b}\right)\right] \\
    & +\sum_{\substack{i+j+2 l \geq 7 \\
    1 \leq i, j \leq 6,1 \leq l \leq 5}} \sqrt{\varepsilon}^{i+j+2 l-7} \frac{\xi^{l}}{l !}\left[\Gamma\left(\widetilde{f}_{i}^{(l)}, f_{j}^{b b}\right)+\Gamma\left(f_{j}^{b b}, \widetilde{f}_{i}^{(l)}\right)\right] \\
   & +\sum_{\substack{i+j+l \geq 7\\
   1\leq i,j\leq 6,1\leq l\leq 5}} \sqrt{\varepsilon}^{i+j+l-7} \frac{\xi^{l}}{l !}\left[\Gamma\left(\widetilde{f}_{i}^{b,(l)}, f_{j}^{b b}\right)+ \Gamma \left(f_{j}^{b b}, \widetilde{f}_{i}^{b,(l)}\right)\right] \\
    &+\frac{\xi^{6}}{6 !} \sum_{j=1}^{6} \sqrt{\varepsilon}^{j-1}\left[\Gamma\left( \sum_{i=1}^{6} \sqrt{\varepsilon}^{i+6} \widetilde{f}_{i}^{(6)}+\sqrt{\varepsilon}^{i} \widetilde{f}_{i}^{b,(6)}, f_{j}^{b b}\right)\right. \\
    &\left.+\Gamma\left(f_{j}^{b b},   \sum_{i=1}^{6} \sqrt{\varepsilon}^{i+6} \widetilde{f}_{i}^{(6)}+\sqrt{\varepsilon}^{i} \widetilde{f}_{i}^{b,(6)}\right)\right] .
    \end{aligned}
\end{multline}
Furthermore, the following initial data are imposed on the remainder equation \eqref{1.5.1}:
\begin{equation}\label{1.5.10}
  f_{R,\varepsilon}(0,x,v)=f^{in}_{R,\varepsilon}(x,v),
\end{equation}
which satisfies the compatibility condition on $\Sigma_-$
\begin{equation*}
  \gamma_-f^{in}_{R,\varepsilon}=\mathcal{K}f^{in}_{R,\varepsilon}.
\end{equation*}

%\subsection{Sketch of ideas.}

%In this paper, we expect to justify the limit process rigorously by the Hilbert expansion method. We employ $L^2-L^\infty$ framework to prove that the remainder $\sqrt{\varepsilon}^5f_{R,\varepsilon}\to 0$ as $\varepsilon\to 0$. The uniform $L^2-L^\infty$ estimates rely on an interplay between $L^2$ and $L^\infty$ estimates for the Boltzmann equation.

%Indeed, the remainder equation has the  following form
%\begin{equation}
 % \partial_t f_{R,\varepsilon}+v\cdot\nabla_xf_{R,\varepsilon}+\frac{1}{\varepsilon}\mathcal{L}f_{R,\varepsilon}=\text{some other terms. }
%\end{equation}
%The coercivity \eqref{label} of the linear Boltzmann operator $\mathcal{L}$ provides a kinetic dissipative rate $\frac{c_0}{\varepsilon}\|(\mathbf{I-P})f_{R,\varepsilon}\|_\nu^2$ in $L^2$ estimates. This causes $\|f_{R,\varepsilon}\|_2$ to be bounded by $\varepsilon^2\|h\|_\infty$. On the other hand, in $L^\infty$ estimates, after transforming $f_{R,\varepsilon} $ to $h$. Then the norm $\varepsilon^{\frac{3}{2}}\|h\|_\infty$ can be bounded by $\varepsilon^\frac{3}{2}\|h(0)\|_\infty+\|f_{R,\varepsilon}\|_2$. However, there is a $\sqrt{\varepsilon}^3-$order disparity between $\sup_{0\leq s\leq t}\|h(s)\\_\infty$ and $\sup_{0\leq s\leq t}\|f_{R,\varepsilon}(s)\|_2$.

\section{Uniform Bounds To the solutions of expansions}
In this section, we focus on the uniform bounds to $F_k,F^b_k$ and $F^{bb}_k(1\leq k\leq 6)$ appears in the Hilbert expansion, which will be applied to control the source term in the remainder equation \eqref{1.5.2}.
\subsection{Some auxiliary systems}

Firstly, we consider the estimates for the interior term $F_k$. We only need to control the fluid variables due to the microscopic parts are totally determined by the macroscopic parts. Let $(\tilde{\rho},\tilde{u},\tilde{\theta})(t,x)$ be the solution to the linear system:
\begin{equation}\label{2.1.1}
  \begin{cases}
    \partial_t\tilde{\rho}+\nabla_x\cdot \tilde{u}=0,\\
    \partial_t u_1+\nabla_x(\rho_1+\theta_1)=f,\\
    \partial_t\tilde{\theta}+\frac{2}{3}\nabla_x\cdot\tilde{u}=g,
  \end{cases}
\end{equation}
with $(t,x)\in(0,\tau)\times \mathbb{R}^3_+$. We impose the following boundary condition
\begin{equation}\label{2.1.2}
  \tilde{u}_3(t,\bar{x},0)=d(t,\bar{x}),\quad\forall (t,\bar{x})\in(0,\tau)\times \mathbb{R}^2,
\end{equation}
and initial condition
\begin{equation}\label{2.1.3}
  (\tilde{\rho},\tilde{u},\tilde{\theta})(0,x)=(\tilde{\rho}_0,\tilde{u}_0,\tilde{\theta}_0).
\end{equation}
By \cite{guo2021hilbert},  the following lemma for the existence of smooth solution for the linear system \eqref{2.1.1} holds:
\begin{lemma}\label{Lm2.1}
  Assume that
  \begin{equation}\label{2.1.4}
    \mathbb{E}_0:=\|(\tilde{\rho}^0,\tilde{u}^0,\tilde{\theta}^0)\|^2_{\mathcal{H}^{k+1}(\mathbb{R}^3_+)}+\sup_{t\in (0,\tau)}\left[\left\|(f,g)(t)\right\|^2_{\mathcal{H}^{k+1}(\mathbb{R}^3_+)}+\|d(t)\|^2_{\mathcal{H}^{k+2}(\mathbb{R}^2)}\right]<+\infty,
  \end{equation}
  with $k\geq 3 $, and the compatibility condition is satisfied for the initial data. Then there exist a unique smooth solution to \eqref{2.1.1} with the boundary condition \eqref{2.1.2} for $t\in[0,\tau],$ such that
  \begin{equation}\label{2.1.5}
    \sup_{t\in[0,\tau]}\|(\tilde{\rho},\tilde{u},\tilde{\theta})(t)\|^2_{\mathcal{H}^k(\mathbb{R}^3_+)}\leq C_\tau \mathbb{E}_0.
  \end{equation}
\end{lemma}

Secondly, the viscous boundary layers $F^b_k$ will be controlled. As shown in Lemma \ref{Lm1.1} and Remark \ref{Rm1.2}, to construct the solution of the viscous boundary layer, we need to solve the following heat equations of $(u,\theta)=(u_1,u_2,\theta)(t,\bar{x},\zeta)$:
\begin{equation}\label{2.1.6}
  \begin{cases}
    \partial_tu_i-\kappa_1\partial^2_\zeta u_i=f_i, \quad (i=1,2),\\
    \partial_t\theta -\frac{2}{5}\kappa_2\partial^2_\zeta\theta=g,
  \end{cases}
\end{equation}
with the Dirichlet boundary condition:
\begin{equation}\label{2.1.7}
  \begin{cases}
    u_i|_{\zeta=0}=b_i(t,\bar{x}),\quad \theta|_{\zeta=0}=a(t,\bar{x}),\\
    \lim_{\zeta\to 0}(u_i,\theta)(t,\bar{x},\zeta)=0\quad (i=1,2).
  \end{cases}
\end{equation}
where $(t,\bar{x},\zeta)\in [0,\tau]\times \mathbb{R}^2\times \mathbb{R}_+$, and $\kappa_1,\kappa_2$ are positive constants defined in \eqref{A.3}. The source terms $(f_1,f_2,g)(t,\bar{x},\zeta)$ are already given. The boundary condition $a(t,\bar{x})$ and  $b_i(t,\bar{x})$ are determined by the fluid quantities of the lower order.  We impose \eqref{2.1.6} with initial data
\begin{equation}
  u(t,\bar{x},\zeta)|_{t=0}=u_0(\bar{x},\zeta),\quad \theta(t,\bar{x},\zeta)|_{t=0}=\theta_0(\bar{x},\zeta),
\end{equation}
which satisfies the corresponding compatibility condition.

Recalling the definitions of $\mathbb{H}^k_l(\mathbb{R}^3_+)$ and $\mathbb{H}^k(\mathbb{R}^2)$ in \eqref{1.3.14} and \eqref{1.3.15}, we have the following lemma.

\begin{lemma}\label{Lm2.2}
  Assume that
  \begin{equation}
    \begin{aligned}
      \mathbb{E}_1:=&\|u_0\|^2_{\mathbb{H}^k_l(\mathbb{R}^3_+)}+\sup_{t\in [0,\tau]}\left(\|b(t)\|_{\mathbb{H}^{k+3}(\mathbb{R}^2)}+\|f(t)\|^2_{\mathbb{H}^{k+1}_l(\mathbb{R}^3_+)}\right)\\
      &+\|\theta_0\|^2_{\mathbb{H}^k_l(\mathbb{R}^3_+)}+\sup_{t\in [0,\tau]}\left(\|a(t)\|_{\mathbb{H}^{k+3}(\mathbb{R}^2)}+\|g(t)\|^2_{\mathbb{H}^{k+1}_l(\mathbb{R}^3_+)}\right)\leq \infty,
    \end{aligned}
  \end{equation}
  with $k\geq 3$, $l\geq 0$.
Then there exists a unique smooth solution $(u,\theta)(t,\bar{x},\zeta)$ to \eqref{2.1.6} over $t\in[0,\tau]$ satisfying
\begin{equation}
  \sup_{t\in[0,\tau]}\left(\|(u,\theta)\|^2_{\mathbb{H}^k_l(\mathbb{R}^3_+)}+\int_{0}^{t}\|\partial_\zeta(u,\theta)\|_{\mathbb{H}^k_l(\mathbb{R}^3_+)}^2ds\right)\leq C_\tau \mathbb{E}_1.
\end{equation}
\end{lemma}
The proof of Lemma \ref{Lm2.2} is similar and much easier as \cite[Lemma 4.1]{guo2021hilbert}. For simplicity, we omit the details of proof here.

We finally introduce the existence results on the solution of the Knudsen-Layer problem. Let $f(t,\bar{x},\xi,v)$ satisfy the following boundary layer problem with $(t,\bar{x},\xi,v)\in[0,\tau]\times \mathbb{R}^2\times \mathbb{R}_+\times \mathbb{R}^3_+:$
\begin{equation}\label{2.1.11}
  \begin{cases}
    v_3\partial_\xi f+\mathcal{L}f=s(t,\bar{x},\xi,v),\\
    \tilde{\mathcal{K}}f=\tilde{\mathcal{K}}q^\infty+g(t,\bar{x},v_3),\\
    \lim_{\xi\to0}f(t,\bar{x},\xi,v)=0,
  \end{cases}
\end{equation}
where
\begin{equation}
  q^\infty=\left[b_1v_1+b_2v_2+c\left(\frac{|v|^2-3}{2}\right)\right]\sqrt{\mu}.
\end{equation}
The Knudsen layer problem \eqref{2.1.11} has been solved in the $L^2(e^{\sigma\xi} d\xi;L^2(|v_3|dv))$ space by Coron-Golse-Sulem \cite{golse1988boundary}. Recently, authors in \cite{huang2022boundary} and \cite{Jiang2024boundary} derived a pointwise estimate for the case accommodation coefficient $\alpha=1$ and $0<\alpha<1$, respectively.

\begin{lemma}\label{Lm2.3}
  Let $0<a<\frac{1}{2}$ and $w_l$ be defined in \eqref{1.5.13} with $l\geq 3$. For each $(t,\bar{x})\in [0,\tau]\times \mathbb{R}^2$, we assume that $s\in \mathcal{N}^\perp$ and the solvability condition $\int_{v_3<0}v_3g\sqrt{\mu}dv=0$ holds, and
  \begin{equation}\label{2.1.12}
    \mathbb{E}_2:=\sup_{t\in [0,\tau]}\sum_{|\alpha|\leq r}\left\{\left\|w_l\mu^{-a}\partial^\alpha_{t,\bar{x}}g(t)\right\|_{L^\infty_{\bar{x},v}\cap L^2_{\bar{x}}L^\infty_v}+\left\|w_l\mu^{-a} e^{\sigma_0 \xi}\partial^\alpha_{t,\bar{x}}s (t)\right\|_{L^\infty_{\bar{x},\xi,v}\cap L^2_{\bar{x}}L^\infty_{\xi,v}}\right\}\leq \infty,
  \end{equation}
for some small $\sigma_0>0$. Then there exists a unique $q^\infty$ with
\begin{equation}\label{2.1.13}
\begin{aligned}
   & q^\infty(v)=[b_1v_1+b_2v_2+c^\infty\left(\frac{|v|^2-3}{2}\right)]\sqrt{\mu},\\
    &|(b_1^\infty,b_2^\infty,c^\infty)|\leq C\mathbb{E}_2,
\end{aligned}
\end{equation}
such that the boundary layer problem \eqref{2.1.11} admits a unique solution $f(t,\bar{x},\xi,v)$ satisfying
\begin{equation}\label{2.1.14}
  \sup_{t\in [0,\tau]}\sum_{|\alpha|\leq r}\left\{\left\|w_l\mu^{-a}\partial^\alpha_{t,\bar{x}}f(t,\cdot,0,\cdot)\right\|_{L^\infty_{\bar{x},v}\cap L^2_{\bar{x}}L^\infty_v}+\left\|w_l\mu^{-a} e^{\sigma_0 \xi}\partial^\alpha_{t,\bar{x}}f(t)\right\|_{L^\infty_{\bar{x},\xi,v}\cap L^2_{\bar{x}}L^\infty_{\xi,v}}\right\}\leq \frac{C}{\sigma_0-\sigma}\mathbb{E}_2,
\end{equation}
for all $\lambda\in (0,\lambda_0)$, where $C>0$ is a constant independent of $(t,\bar{x})$.
\end{lemma}

We note that there is no normal derivatives estimates for the boundary value problem. However, it is sufficient to complete the prove of the main theorem with the tangential and time derivatives estimates in Lemma \ref{Lm2.3}.

Since the source term $s\in \mathcal{N}^{\perp}$ is demanded in Lemma \ref{Lm2.3}, but $S^{bb}_k\notin \mathcal{N}^\perp$ in general. We can decompose $S^{bb}_k$ into microscopic part  and macroscopic part, as shown in \eqref{Sbb-1} and \eqref{Sbb-2}. We also emphasize that when solving $f^{bb}_k$, the source terms $S^{bb}_{k,1}$ and $S^{bb}_{k,2}$ are known, as $f_i,f_i^b$ and $f^{bb}_j(i\leq k,j\leq k-1)$ have already been determined.
 We then apply the same arguments as \cite[Lemma 2.9]{guo2021hilbert} to cancel $S^{bb}_{k,1}$ in \eqref{1.4.3.11}. Thus, we can always use Lemma \ref{Lm2.3} in the Knudsen layer problem \eqref{2.1.11}.

\begin{lemma}
 We assume that
\begin{equation}
   S^{bb}_{k,1}=\{a_k+b_k\cdot v+c_k|v|^2\}\sqrt{\mu}
\end{equation}
satisfy
\begin{equation*}
  \lim_{\xi\to\infty}e^{\eta\xi}|(a_k,b_k,c_k)(t,\bar{x},\xi)|=0
\end{equation*}
for some positive constant $\eta>0$. Then there exists a function
\begin{equation*}
  f_{k, 1}^{b b}=\left\{\Psi_{k} v_{3}+\Phi_{k, 1} v_{3}v_1+\Phi_{k, 2} v_{3}v_2+\Phi_{k, 3}+\Theta_{k} v_{3}|v|^2\right\} \sqrt{\mu}
\end{equation*}
  such that  $v_{3} \partial_{\xi} f_{k, 1}^{b b}-S_{k, 1}^{b b} \in \mathcal{N} ^{\perp}$, where

\begin{equation}
    \begin{array}{l}
    \Psi_{k}(t, \bar{x}, \xi)=-\int_{\xi}^{+\infty}\left(2 a_{k}+3 c_{k}\right)(t, \bar{x}, s) d  s, \\
    \Phi_{k, i}(t, \bar{x}, \xi)=-\int_{\xi}^{+\infty}   b_{k, i}(t, \bar{x}, s) d  s,\quad  i=1,2, \\
    \Phi_{k, 3}(t, \bar{x}, \xi)=-\int_{\xi}^{+\infty} b_{k, 3}(t, \bar{x}, s) d  s ,\\
    \Theta_{k}(t, \bar{x}, \xi)=\frac{1}{5 } \int_{\xi}^{+\infty} a_{k}(t, \bar{x}, s) d  s.
    \end{array}
\end{equation}
  Moreover, there holds
\begin{equation}
    \begin{array}{l}
    \left|v_{3} \partial_{\xi} f_{k, 1}^{b b}-S_{k, 1}^{b b}\right| \leq C\left|\left(a_{k}, b_{k}, c_{k}\right)(t, \bar{x}, \xi)\right|\langle v\rangle^{4} \sqrt{\mu}, \\
    \left|f_{k, 1}^{b b}(t, \bar{x}, \xi, v)\right| \leq C\langle v\rangle^{3} \sqrt{\mu} \int_{\xi}^{\infty}\left|\left(a_{k}, b_{k}, c_{k}\right)\right| \rightarrow 0 \text { as } \xi \rightarrow \infty .
    \end{array}
\end{equation}
\end{lemma}

   It is known from Subsection \ref{se1.4.3} that $S^{bb}_{2,1}=S^{bb}_{3,1}=0$. Denote $f^{bb}_{k,2}=f^{bb}_k-f^{bb}_{k,1}$. We thereby have
\begin{equation}
    v_{3} \partial_{\xi} f_{k, 2}^{b b}+\mathcal{L} f_{k, 2}^{b b}=S_{k, 2}^{b b}-\left(v_{3} \partial_{\xi} f_{k, 1}^{b b}+\mathcal{L}f^{bb}_{k,1}-S_{k, 1}^{b b}\right) \in \mathcal{N}  ^{\perp},\quad \lim _{\xi \rightarrow \infty} f_{k, 2}^{b b}(t, \bar{x}, \xi, v)=0.
\end{equation}
Therefore, one can solve $f^{bb}_k$ by solving $f^{bb}_{k,2}$.

\subsection{Uniform bounds for $F_k,F^b_k,F_k^{bb}(1\leq k\leq 6)$}

Based on Lemma \ref{Lm2.1}, Lemma \ref{Lm2.2} and Lemma \ref{Lm2.3}, following the analogously arguments in \cite[Proposition 5.1]{guo2021hilbert}, we obtain the following uniform bounds of the expansions $F_k,F^b_k,F^{bb}_k$ $(1\leq k\leq 6)$. We omit the details of proof here for brevity.

\begin{proposition} \label{Pro2.4}
Let $0<a<\frac{1}{2}$. There are $s_k,s^b_k,s^{bb}_k\in \mathbb{N}_+,p_k,p_k^b,p^{bb}_k\in\mathbb{R}_+$ for $1\leq k\leq 6$ satisfying
\begin{equation*}\label{2.2.1}
\begin{aligned}
  s^{bb}_1=0, s_1=s_1^b\gg 1,\\
    s_1> s_k>s^b_k>s^{bb}_k\geq s_{k+1}>s^b_{k+1}>s^{bb}_{k+1}\geq\dots \gg 1\quad   (2\leq k\leq 5),\\
    p_k\gg p^b_k\gg p^{bb}_k\gg p_{k+1}\gg 1\quad (1\leq k\leq 5).
\end{aligned}
\end{equation*}
and $l^k_j=l^b_k+2(s^b_k-j)$ $(0\leq j \leq s^b_k)$ with
\begin{equation}\label{2.2.2}
  l^6_j\geq 10,\,l^k_j\geq 2l^{k+1}_j+28\,(1\leq k\leq 5),
\end{equation}
such that if the initial data $(\rho^{in}_k,u^{in}_k,\theta^{in}_k)(1\leq k\leq 6)$ in \eqref{1.4.1.13} and $((u^{b,in}_k,\theta^{b,in}_k))(1\leq k\leq 6)$ in \eqref{1.4.2.31} satisfy \eqref{1.5.14}, then there exist solutions $f_k,f_k^b,f^{bb}_k$ constructed in subsection \ref{Se1.4} over time interval $t\in [0,\tau]$, respectively, satisfying
\begin{equation}\label{2.2.3}
\begin{aligned}
    \sup_{t\in [0,\tau]}\sum_{k=1}^{6}&\left\{\sum_{\gamma+|\beta|\leq s_k}\left\|w_{p_k}\mu^{-a}\partial_t^\gamma\partial^\beta_xf_k(t)\right\|_{L^2_xL^\infty_v}+\sum_{j=0}^{s^b_k}\sum_{2\gamma+|\bar{\beta}=j}\left\|w_{p^b_k}\mu^{-a}\partial_t^\gamma\partial^{\bar{\beta}}_{\bar{x}}f^b_k(t)\right\|_{L^2_{l^k_j}L^\infty_v}\right.\\
    &\left.+\sum_{\gamma+|\bar{\beta}|\leq s^{bb}_k}\left\|e^{\frac{\xi}{2^{k-1}}}w_{p^{bb}_k}\mu^{-a}\partial_t^\gamma\partial^{\bar{\beta}}_{\bar{x}}f^{bb}_k(t)\right\|_{L^\infty_{\bar{x},\xi,v}\cap L^2_{\bar{x}}L^\infty_{\xi,v}}
    \right\}\\
    &\leq C(\tau,\mathcal{E}^{in}),
\end{aligned}
\end{equation}
where the notation $f^{bb}_1=0$ is employed.
\end{proposition}
We emphasize that the viscous boundary layers $F^b_k$ decay algebraically associated with $\zeta$, and the Knudsen boundary layers $F^{bb}_k$ decay exponentially associated with $\xi>0$. These suffice to dominate the source terms in \eqref{1.5.1} while deriving the uniform estimates for the remainder $f_{R,\varepsilon}$.
%\begin{corollary}\label{Co2.5}
 % Under the same assumptions as in Proposition \ref{Pro2.4}, we have the following estimates:
 %\begin{equation}
  %\begin{aligned}
   % \|S_2\|^2_2\leq \frac{C}{\sqrt{\varepsilon}}\|f_{R,\varepsilon}\|^2_\nu \quad \|w_lS_2\|_\infty\leq \frac{C}{\sqrt{\varepsilon}}\|h_{R,\varepsilon}\|_{\infty},\\
%\|S_i\|_{L^2}\leq C,\quad \|w_lS_i\|_{\infty}\leq C \quad (i=3,4,5),
 % \end{aligned}
 %\end{equation}
 %where $S_i$ $(i=2,3,4,5)$ are defined in \eqref{1.5.6}-\eqref{1.5.9} and $C=C(\tau,\mathcal{E}^{in})> 0.$
%\end{corollary}

\section{Uniform Bounds for Remainder $f_{R,\varepsilon}$}
This section is dedicated to the $L^2$-$L^\infty$ estimates for the solutions of the remainder equation \eqref{1.5.2}. This is sufficient to obtain the uniform estimates of $f_{R,\varepsilon}$.

\begin{lemma}\label{Lm3.1} Under the same assumptions as Proposition \ref{Pro2.4}. Let $f_{R,\varepsilon}$ be the solution of \eqref{1.5.2} and $h_{R,\varepsilon}$ be defined as in \eqref{1.5.12} with $l\geq 7$. For any given $\tau>0$ there exist constants $\varepsilon'_0>0$ and $C=C(\tau,\mathcal{E}^{in})> 0$ such that for all $0<\varepsilon<\varepsilon'_0$ and $t\in [0,\tau]$,
\begin{multline}\label{3.1}
      \frac{d}{dt}\|f_{R,\varepsilon}\|_2^2+\frac{c_0}{2\varepsilon}\left\|(\mathbf{I-P})f_{R,\varepsilon}\right\|_\nu^2+\frac{(2-\alpha)\alpha}{2}\int_{\gamma_+}(v\cdot n)\left|(I-P_\gamma)f_{R,\epsilon}\right|^2d\bar{x}dv\\
      \leq C\left(\varepsilon^4\|h_{R,\varepsilon}\|^2_\infty+1\right)\|f_{R,\varepsilon}\|_2^2+C,
\end{multline}
  where the positive constant $c_0$ is given in \eqref{1.3.6} and the boundary operator $P_\gamma$ is defined in \eqref{3.1.4}
\end{lemma}

\begin{lemma}\label{Lm3.2}
  Under the same assumptions as Lemma \ref{Lm3.1}, there are constants $\varepsilon''_0$ and $C=C(\tau,\mathcal{E}^{in})$ such that for all $0< \varepsilon<\varepsilon''_0$ and $t\in[0,\tau]$,
  \begin{equation}\label{3.2}
    \sup_{0\leq s\leq t}\left\|\sqrt{\varepsilon}^3h_{R,\varepsilon}(s)\right\|_\infty\leq C\left(\left\|\sqrt{\varepsilon}^3h_{R,\varepsilon}(0)\right\|_\infty+\sup_{0\leq s\leq \tau}\left\|f_{R,\varepsilon}(s)\right\|_2+\sqrt{\varepsilon}^5\right).
  \end{equation}
\end{lemma}
\begin{proof}
  [Proof of Theorem \ref{Thm1.4}.] With Lemma \ref{Lm3.1} and Lemma \ref{Lm3.2}, the rest of the proof is the same as \cite{guo2010acoustic,jiang2021compressibleb}. For simplicity, we omit the details here. Therefore we conclude the proof of Theorem \ref{Thm1.4}.
\end{proof}

\subsection{$L^2$ estimates: Proof of Lemma \ref{Lm3.1}.}
Multiplying \eqref{1.5.2} by $f_{R,\varepsilon}$ and integrating respect to $(x,v)$ over $\mathbb{R}^3_+\times \mathbb{R}^3$, one obtains
\begin{equation}\label{3.1.1}
  \frac{1}{2}\frac{d}{dt}\|f_{R,\varepsilon}\|_2^2+\frac{1}{\varepsilon}\left<\mathcal{L}f_{R,\varepsilon},f_{R,\varepsilon}\right>+\frac{1}{2}\int_{\Sigma}(v\cdot n)|f_{R,\varepsilon}|^2d\bar{x}dv=\left<S,f_{R,\varepsilon}\right>.
\end{equation}
Notice that
\begin{equation}\label{3.1.2}
  \left<\mathcal{L}f_{R,\varepsilon},f_{R,\varepsilon}\right>\geq c_0\left\|(\mathbf{I-P})f_{R,\varepsilon}\right\|_\nu^2.
\end{equation}
Using the boundary condition \eqref{1.5.3}, the boundary term on the LHS can be estimated as
\begin{equation}\label{3.1.3}
\begin{aligned}
  &\frac{1}{2} \int_{\Sigma}(v\cdot n)|f_{R,\varepsilon}|^2d\bar{x}dv\\
  & =\frac{1}{2}\int_{\gamma_-}(v\cdot n)\left[(1-\alpha)L\gamma_+f_{R,\varepsilon}+\alpha P_\gamma f_{R,\varepsilon}\right]^2d\bar{x}dv+\frac{1}{2}\int_{\gamma_+}(v\cdot n)f_{R,\varepsilon}^2d\bar{x}dv\\
  & =\frac{(2-\alpha)\alpha}{2}\int_{\gamma_+}(v\cdot n)f_{R,\epsilon}^2d\bar{x}dv-\frac{(2-\alpha)\alpha}{2}\sqrt{2\pi}\int_{\partial\mathbb{R}^3_+}\left(\int_{v\cdot n>0}(v\cdot n)f_{R,\epsilon}\sqrt{\mu(v)}dv\right)^2d\bar{x}\\
  & =\frac{(2-\alpha)\alpha}{2}\int_{\gamma_+}(v\cdot n)\left|(I-P_\gamma)f_{R,\epsilon}\right|^2d\bar{x}dv,
\end{aligned}
\end{equation}
where we used the notation
\begin{equation}\label{3.1.4}
  P_\gamma f(v)=\sqrt{2\pi \mu(v)}\int_{v'\cdot n>0}(v'\cdot n)f(v')\sqrt{\mu(v')}dv'.
\end{equation}
For the source term on the RHS of \eqref{1.5.2}, we split
\begin{equation}\label{3.1.5}
  \left<S,f_{R,\varepsilon}\right>= \left<S,\mathbf{P}f_{R,\varepsilon}\right>+ \left<S,(\mathbf{I-P})f_{R,\varepsilon}  \right>.
\end{equation}
Notice that $\Gamma(f,g)\in \mathcal{N}^\perp$, one has
\begin{equation}\label{3.1.6}
\begin{aligned}
   \left<S_1,f_{R,\varepsilon}\right>&= \sqrt{\varepsilon}^3\int_{\mathbb{R}^3_+\times \mathbb{R}^3} \Gamma(f_{R,\varepsilon},f_{R,\varepsilon})f_{R,\varepsilon}dxdv\\
  &=\sqrt{\varepsilon}^3\int_{\mathbb{R}^3_+\times \mathbb{R}^3} \Gamma(f_{R,\varepsilon},f_{R,\varepsilon})(\mathbf{I-P})f_{R,\varepsilon}dxdv\\
  &\leq \sqrt{\varepsilon}^3\|(\mathbf{I-P})f_{R,\varepsilon}\|_\nu\|h_{R,\varepsilon}\|_{L^\infty}\|f_{R,\varepsilon}\|_{L^2}\\
  &\leq \frac{\delta}{\varepsilon} \|(\mathbf{I-P})f_{R,\varepsilon}\|^2_\nu+\varepsilon^4C_\delta\|h_{R,\varepsilon}\|^2_{L^\infty}\|f_{R,\varepsilon}\|^2_{L^2}.
\end{aligned}
\end{equation}
It's known from Proposition \ref{Pro2.4} and \eqref{1.5.6}-\eqref{1.5.9} that
\begin{equation}\label{3.1.7}
  \begin{aligned}
     \left<S_2,f_{R,\varepsilon}\right>&=\int_{\mathbb{R}^3_+\times \mathbb{R}^3}S_2(\mathbf{I-P})f_{R,\varepsilon}dxdv\\
    &\leq  \frac{C}{\sqrt{\varepsilon}}\left\|(\mathbf{I-P})f_{R,\varepsilon}\right\|_\nu\|f_{R,\varepsilon}\|_\nu\\
    &\leq \frac{\delta}{\varepsilon} \left\|(\mathbf{I-P})f_{R,\varepsilon}\right\|^2_\nu+ C_\delta\|f_{R,\varepsilon}\|_\nu^2\\
    &\leq (\delta+C_\delta \varepsilon)\frac{1}{\varepsilon}\|(\mathbf{I-P})f_{R,\varepsilon}\|_\nu^2+C_\delta\|f_{R,\varepsilon}\|^2_2,
  \end{aligned}
  \end{equation}
and
  \begin{equation}\label{3.1.8}
       \int_{\mathbb{R}^3_+\times \mathbb{R}^3} S_if_{R,\varepsilon}dxdv\leq C \| f_{R,\varepsilon}\|_2,\quad (i=3,4,5).
    \end{equation}
Collecting \eqref{3.1.5}-\eqref{3.1.8}, one obtains
\begin{equation}\label{3.1.9}
  \left<S,f_{R,\varepsilon}\right>\leq (\frac{2\delta}{\varepsilon}+C\varepsilon)\|(\mathbf{I-P})f_{R,\varepsilon}\|_\nu^2+\varepsilon^4C_\delta\|h_{R,\varepsilon}\|^2_\infty\|f_{R,\varepsilon}\|_2^2+C_\delta\|f_{R,\varepsilon}\|^2_2+C.
\end{equation}
Hence \eqref{3.1} follows from \eqref{3.1.1}-\eqref{3.1.3} and \eqref{3.1.9} by choosing $\varepsilon>0$ and $\delta>0$ suitable small. This completes the proof of Lemma \ref{Lm3.1}.

\subsection{$L^\infty$ estimates: Proof of Lemma \ref{Lm3.2}.}
In this subsection, we deduce the $L^\infty$ estimates for the solutions of the remainder equation \eqref{1.5.2}.  As in \cite{guo2010decay} or \cite{guo2021hilbert}, we introduce the following backward characteristics. Given $(t,x,v)$, we define $[X(s),V(s)]$ satisfy
\begin{equation*}
  \begin{cases}
    \frac{d}{ds}X(s)=V(s),\quad \frac{d}{dx}V(s)=0,\\
    [X(t;t,x,v),V(t;t,x,v)]=[x,v].
  \end{cases}
\end{equation*}
Then $[X(s;t,x,v),V(s;t,x,v)]=[x-(t-s)v,v]=[X(s),V(s)]$. For $(x,v)\in \mathbb{R}^3_+\times \mathbb{R}^3$, the backward exist time $t_b(x,v)>0$ is defined to be the first moment at which the backward characteristics line $[X(s;0,x,v);V(s;0,x,v)]$ emerges from $\partial\mathbb{R}^3_+$. Hence, it holds that
$$t_b(x,v)=\inf\{t>0:x-tv\notin \mathbb{R}^3_+\},$$
which indicates that $x-t_b v\in \partial\mathbb{R}^3_+$. We also define $x_b(x,v)=x-t_bv$. Note that we use $t_b(x,v)$ whenever it is well-defined.
For half space problem, the backward trajectory does not hit the boundary for the case $v_3<0$, and the backward trajectory will hit the boundary once for the case $v_3>0$. By the Maxwell reflection boundary condition, the backward characteristics can be written as
\begin{equation*}
\begin{aligned}
    &X(s;t,x,v)=\mathbbm{1}_{[t_1,t)}(s)\{x-(t-s)v\}+\mathbbm{1}_{(-\infty,t_1)}(s)\{x-[(1-\alpha)R_{x_1}v+\alpha v'](t-s)\},\\
    &V(s;t,x,v)=\mathbbm{1}_{[t_1,t)}(s)v+\mathbbm{1}_{(-\infty,t_1)}(s)[(1-\alpha)R_{x_1}v+\alpha v'],
\end{aligned}
\end{equation*}
where $v'\in \mathcal{V}:=\{v':v'\cdot n>0\}$ and
$$(t_1,x_1)=(t-t_b(x,v),x_b(x,v)).$$

\begin{lemma}\label{Lm3.3}\cite[Lemma 3]{guo2010decay}
Recall \eqref{1.3.7}. We have
  \begin{equation}\label{1.3.8}
   |k(v, v')| \leq C\{|v-v'|+|v-v'|^{-1}\} e^{-\frac{|v-v'|^{2}}{8}} e^{-\frac{\left||v|^{2}-|v'|^{2}\right|^{2}}{8|v-v'|^{2}}},
  \end{equation}
   Let $0\leq a<\frac{1}{4}$ and $l\geq 0$. Then for $\eta>0$ sufficiently small.
  \begin{equation}
    \int_{\mathbb{R}^3}e^{\eta|v-v'|^2}|k(v,v')|\frac{w_l(v) e^{a|v|^2}}{w_l(v') e^{a|v'|^2}}dv'\leq \frac{C}{1+|v|},
  \end{equation}
  for some $ C>0$.
\end{lemma}
%Following (2.1), it is direct to have\begin{equation*}\int_{\mathbb{R}^{3}}\left|k(v, v') \cdot \frac{(1+|v|)^%{\alpha}}{(1+|v'|)^{\alpha}}\right| d v'\leq C_{\alpha}(1+|v|)^{-1}\end{equation*}

We denote the weighted non-local operator
\begin{equation}\label{3.2.1}
  K_wg=wK(\frac{g}{w})=\int_{\mathbb{R}_3}k_w(v,v')g(v')dv',
\end{equation}
where
\begin{equation}
  k_w(v,v'):=k(v,v')\frac{w_l(v)}{w_l(v')}.
\end{equation}
Based on the remainder equation \eqref{1.5.2} and the definition \eqref{1.5.12}, we deduce that $h_{R,\varepsilon}$ satisfies
\begin{equation}\label{3.2.2}
  \begin{cases}
        \partial_t h_{R,\varepsilon}+v\cdot\nabla_xh_{R,\varepsilon}+\frac{\nu}{\varepsilon}h_{R,\varepsilon}-\frac{1}{\varepsilon}K_wh_{R,\varepsilon}=w_lS,\\
        h_{R,\varepsilon}(t,x,v)=\int_{v'\cdot n>0}\mathcal{R}(v,v')h_{R,\varepsilon}(t,x,v')dv' \quad (x,v)\in \Sigma_-,
  \end{cases}
  \end{equation}
where we used the notation
\begin{equation}\label{3.2.3}
\mathcal{R}(v,v')=(1-\alpha)\delta(v'-(v-2v\cdot n)n)+\alpha\sqrt{2\pi}\sqrt{\mu(v)\mu(v')}|v'\cdot n|\frac{w_l(v)}{w_l(v')}.
\end{equation}
We can integrate $\eqref{3.2.2}_1$ along the backward characteristics to obtain
\begin{equation}\label{3.2.4}
  \begin{aligned}
      h_{R,\varepsilon}(t,x,v)=&\mathbbm{1}_{t_1\leq 0}\Bigg\{e^{-\frac{\nu t}{\varepsilon}}h_{R,\varepsilon}(0,x-tv,v) \\
      &+ \int_{0}^{t} e^{-\frac{\nu(t-s)}{\varepsilon}} \left[\frac{1}{\varepsilon}K_w h_{R,\varepsilon} + w_l S\right](s,x+(s-t)v,v) ds \Bigg\} \\
      &+ \mathbbm{1}_{t_1>0}\Bigg\{e^{\frac{\nu(t_1-t)}{\varepsilon}} h_{R,\varepsilon}(t_1,x_1,v) \\
      &+ \int_{t_1}^{t} e^{-\frac{\nu(t-s)}{\varepsilon}} \left[\frac{1}{\varepsilon}K_w h_{R,\varepsilon} + w_l S\right](s,x+(s-t)v,v) ds \Bigg\}.
  \end{aligned}
  \end{equation}
Using the boundary condition $\eqref{3.2.2}_2$ and \eqref{3.2.4} again to obtain
\begin{equation}\label{3.2.5}
  h_{R,\varepsilon}(t,x,v)=\sum_{i=1}^{5}\mathcal{I}_{i},
\end{equation}
with
\begin{equation}
 \begin{aligned}
   \mathcal{I}_1 =&\mathbbm{1}_{t_1\leq 0}e^{-\frac{\nu t}{\varepsilon}}h_{R,\varepsilon}(0,x-tv,v)+\mathbbm{1}_{t_1>0}e^{-\frac{\nu t}{\varepsilon}}\int_{v_1\cdot n>0}\mathcal{R}(v,v_1)h_{R,\varepsilon}(0,x_1-t_1v_1,v_1)dv_1, \\
   \mathcal{I}_2 =&\int_{\max\{0,t_1\}}^{t}e^{\frac{\nu(s-t)}{\varepsilon}}w_lS(s,x+(s-t)v,v)ds\\
    &+\mathbbm{1}_{t_1>0}\int_{v_1\cdot n>0}\int_{0}^{t_1}\mathcal{R}(v,v_1)e^{\frac{\nu(s-t)}{\varepsilon}}w_lS(s,x_1+(s-t_1)v_1,v_1)dsdv_1, \\
   \mathcal{I}_3 =&\int_{\max\{0,t_1\}}^{t}e^{\frac{\nu(s-t)}{\varepsilon}}\frac{1}{\varepsilon}K_wh_{R,\varepsilon}(s,x+(s-t)v,v)ds, \\
 \end{aligned}
\end{equation}
and
\begin{equation}\label{3.2.9}
  \mathcal{I}_4=\mathbbm{1}_{t_1>0}\int_{v_1\cdot n>0}\int_{0}^{t_1}\mathcal{R}(v,v_1)e^{\frac{\nu(s-t)}{\varepsilon}}\frac{1}{\varepsilon}K_wh(s,x_1+(s-t_1)v_1,v_1)dsdv_1.
\end{equation}
It is direct to see that the term $\mathcal{I}_1$ can be bounded by
\begin{equation}\label{3.2.10}
  |\mathcal{I}_1|\leq C\|h_{R,\varepsilon}(0)\|_\infty.
\end{equation}
We note that (see \cite[Lemma 5]{guo2010decay} for instance)
\begin{equation*}
  |w_l(v)\Gamma(f_{R,\varepsilon},f_{R,\varepsilon})|\leq C\nu\|h_{R,\varepsilon}\|_\infty^2,
\end{equation*}
which together with Proposition \ref{Pro2.4}, one has
\begin{equation}\label{3.2.11}
  |\mathcal{I}_2|\leq  C(\sqrt{\varepsilon}^5 \sup_{0\leq s\leq t}\|h_{R,\varepsilon}(s)\|_\infty^2+\sqrt{\varepsilon}\sup_{0\leq s\leq t}\|h_{R,\varepsilon}(s)\|_\infty+\varepsilon).
\end{equation}
In what follows we only give an explicit computation for the non-local term $\mathcal{I}_4$. The estimates for $\mathcal{I}_3$ are similar and much easier and will be omitted for brevity.

Estimates for $\mathcal{I}_4$. Let us now denote $(t_0',x_0',v_0')=(s,X_{cl}(s;t_1,x_1,v_1),v')$, and define a new back-time cycle as
$$(t_1',x_1',v_1')=(t'_1-t_b(x'_0,v_0'),x_b(x'_0,v'_0),v'_1),$$
for $v_1'\in \mathcal{V'}=\{v_1':v_1'\cdot n>0\}$. We then employ \eqref{3.2.4} again to obtain

\begin{equation}\label{3.2.12}
\begin{aligned}
    \mathcal{I}_4&=\mathbbm{1}_{t_1>0}\int_{v_1\cdot n>0}\int_{0}^{t_1}\mathcal{R}(v,v_1)e^{\frac{\nu(s-t)}{\varepsilon}}\frac{1}{\varepsilon}\int_{\mathbb{R}^3}k_w(v_1,v')h_{R,\varepsilon}(s,x_1+(s-t_1)v_1,v')dsdv'dv_1\\
    &=\sum_{i=1}^{4}\mathcal{J}_i,
\end{aligned}
\end{equation}
with
\begin{equation}\label{3.2.13}
  \begin{aligned}
      \mathcal{J}_1 = & \mathbbm{1}_{t_1>0} \int_{v_1\cdot n>0} \int_{0}^{t_1} \mathcal{R}(v,v_1) e^{\frac{\nu(s-t)}{\varepsilon}} \frac{1}{\varepsilon} \int_{\mathbb{R}^3} k_w(v_1,v')\times \left\{ \mathbbm{1}_{t_1'\leq 0} e^{-\frac{\nu t_0'}{\varepsilon}} h_{R,\varepsilon}(0,x'_0-t'_0v'_0,v'_0) \right. \\
      &\qquad + \left. \mathbbm{1}_{t'_1>0} e^{-\frac{\nu(t_1'-t_0')}{\varepsilon}} \int_{v_1'\cdot n>0} \mathcal{R}(v_0',v_1') h_{R,\varepsilon}(0,x_1'-t_1'v_1',v_1') dv_1' \right\} ds dv' dv_1,
  \end{aligned}
  \end{equation}
\begin{multline}\label{3.2.14}
  \begin{aligned}
      \mathcal{J}_2 = & \mathbbm{1}_{t_1>0}\int_{v_1\cdot n>0}\int_{0}^{t_1}\mathcal{R}(v,v_1)e^{\frac{\nu(s-t)}{\varepsilon}}\frac{1}{\varepsilon}\int_{\mathbb{R}^3}k_w(v_1,v') \\
      & \times \int_{\max\{0,t_1'\}}^{t_0'} e^\frac{\nu(s'-t_0')}{\varepsilon} w_l S\left(s',x'_0+(s'-t'_0)v_0',v'_0\right) ds' ds dv' dv_1 \\
      & + \mathbbm{1}_{t_1>0,t_1'>0}\int_{v_1\cdot n>0}\int_{0}^{t_1}\mathcal{R}(v,v_1)e^{\frac{\nu(s-t)}{\varepsilon}}\frac{1}{\varepsilon}\int_{\mathbb{R}^3}k_w(v_1,v') \\
      & \times \int_{v_1'\cdot n>0}\int_{0}^{t_1'}\mathcal{R}(v_0',v_1')e^{\frac{\nu(s'-t_0')}{\varepsilon}} w_l S(s',x_1'+(s'-t'_1)v_1',v_1') ds' ds dv' dv_1'dv_1,
  \end{aligned}
  \end{multline}
\begin{multline}\label{3.2.15}
      \mathcal{J}_3=\mathbbm{1}_{t_1>0}\int_{v_1\cdot n>0}\int_{0}^{t_1}\mathcal{R}(v,v_1)e^{\frac{\nu(s-t)}{\varepsilon}}\frac{1}{\varepsilon^2}\int_{\mathbb{R}^3}k_w(v_1,v')\int_{\max\{0,t_1'\}}^{t_0'}e^\frac{\nu(s'-t_0')}
      {\varepsilon}\times\\
      K_wh_{R,\varepsilon}(s',x'_0+(s'-t'_0)v_0',v'_0)ds'dsdv'dv_1,
\end{multline}
 
and
  \begin{multline}\label{3.2.16}
    \mathcal{J}_4=\mathbbm{1}_{t_1>0,t_1'>0}\int_{v_1\cdot n>0}\int_{0}^{t_1}\mathcal{R}(v,v_1)e^{\frac{\nu(s-t)}{\varepsilon}}\frac{1}{\varepsilon^2}\int_{\mathbb{R}^3}k_w(v_1,v') \int_{v_1'\cdot n>0}\int_{0}^{t_1'}e^\frac{\nu(s'-t_0')}{\varepsilon}\times\\
    \mathcal{R}(v_0',v_1')K_wh(s',x_1'+(s'-t'_1)v_1',v_1')ds'dsdv'dv_1'dv_1.
\end{multline}
Obviously, one has
\begin{equation}\label{3.2.17}
  |\mathcal{J}_1|\leq \|h_{R,\varepsilon}(0)\|_\infty.
\end{equation}
From Proposition \ref{Pro2.4}, one obtains
\begin{equation}\label{3.2.18}
  |\mathcal{J}_2|\leq C(\sqrt{\varepsilon}^5 \sup_{0\leq s\leq t}\|h_{R,\varepsilon}(s)\|_\infty^2+\sqrt{\varepsilon}\sup_{0\leq s\leq t}\|h_{R,\varepsilon}(s)\|_\infty+\varepsilon).
\end{equation}
Next, we only compute $\mathcal{J}_4$ because the estimates for $\mathcal{J}_3$ are similar and easier. We divide our computations into the following three cases:

Case 1. $|v_1|\geq N_0$ or $ |v_1'|\geq N_0$ with $N_0$ suitable large. It follows from lemma \ref{Lm3.3} that
\begin{equation}\label{3.2.20}
  \int_{\mathbb{R}^3}|k_w(v_1,v')|dv'\leq \frac{C}{N_0},
\end{equation}
or
\begin{equation}\label{3.2.21}
  \int_{\mathbb{R}^3}|k_w(v_1',v'')|dv''\leq \frac{C}{N_0},
\end{equation}
which together with \eqref{3.2.3}, yields that
\begin{equation}\label{3.2.22}
  |\mathcal{J}_4|\leq \frac{C}{N_0}\sup_{0\leq s\leq t}\|h_{R,\varepsilon}(s)\|_\infty.
\end{equation}

Case 2. $|v_1|\leq N_0$ and $|v'|\geq 2N_0$ or $|v_1'|\leq N_0$ and $|v''|\geq 2N_0$. Notice that we have either $|v'-v_1|\geq N_0$ or $|v_1'-v''|\geq N_0$, and thus either of the following holds
\begin{equation}\label{3.2.23}
 | k_w(v_1,v')|\leq C e^{-\frac{\eta N_0^2}{8}}| k_w(v_1,v')|e^{\frac{\eta|v_1-v'|^2}{8}},
\end{equation}
or
\begin{equation}\label{3.2.24}
  |k_w(v_1',v'')|\leq C e^{-\frac{\eta N_0^2}{8}} |k_w(v_1',v'')|e^{\frac{\eta|v_1'-v''|^2}{8}},
\end{equation}
for some small $\eta$, and therefore we have
\begin{equation}\label{3.2.25}
  |\mathcal{J}_4|\leq C_\eta e^{-\frac{\eta}{8}N_0^2}\sup_{0\leq s\leq t}\|h(s)\|_\infty.
\end{equation}

Case 3a. $|v_1|\leq N_0,|v'|\leq 2N_0$ and $|v_1'|\leq N_0,|v''|\leq 2N_0$, $s'\geq t_1'-\varepsilon\kappa_*$,

\begin{equation}\label{3.2.26}
  |\mathcal{J}_4|\leq C_{N_0}\kappa_*\sup_{0\leq s\leq t}\|h_{R,\varepsilon}(s)\|_\infty.
\end{equation}

Case 3b. $|v_1|\leq N_0,|v'|\leq 2N_0$ and $|v_1'|\leq N_0,|v''|\leq 2N_0$, $s'\leq t_1'-\varepsilon\kappa_*$,
\begin{equation}\label{3.2.27}
  \begin{aligned}
      \mathcal{J}_4 = & \frac{1}{\varepsilon^2} \mathbbm{1}_{t_1>0,t_1'>0} \int_{{v_1\cdot n>0,|v_1|\leq N_0}} \int_{v_1'\cdot n>0,|v_1'|\leq N_0} \int_{0}^{t_1} \int_{0}^{t_1'-\varepsilon k_*} e^{\frac{\nu(s-t)}{\varepsilon}} e^{\frac{\nu(s'-t'_0)}{\varepsilon}} \\
      & \times \int_{\substack{|v'|\leq 2N_0}} \int_{\substack{|v''|\leq 2N_0}} k_w(v_1,v') k_w(v_1',v'') \mathcal{R}(v,v_1) \mathcal{R}(v_0',v_1') \\
      & \times h_{R,\varepsilon}(s',x_1'+(s'-t_1')v_1',v'') \, ds' \, ds \, dv' \, dv_1' \, dv_1 \, dv''.
  \end{aligned}
  \end{equation}
  
Based on the Lemma \ref{Lm3.3}, $k_w(v_1,v'$) has a possible integrable singularity of $\frac{1}{|v_1-v'|}$. We choose a number $m(N_0)$ to define
\begin{equation*}
  k_{m,w}(p,v')=\mathbbm{1}_{|p-v'|\geq \frac{1}{m},|v'|\leq m}k_w(p,v')
\end{equation*}
such that
\begin{equation}\label{3.2.28}
  \sup_{|p|\leq 3N_0}\int_{|v'|\leq 3N_0}|k_{m,w}(p,v')-k_w(p,v')|dv'\leq \frac{1}{N_0}.
\end{equation}
We split
\begin{equation*}
\begin{aligned}
  k_w(v_1,v')k_w(v_1',v'')  =&\{k_w(v_1,v')-k_{m,w}(v_1,v')\}k_{w}(v_1',v'')\\
  &+\{k_w(v_1',v'')-k_{m,w}(v_1',v'')\}k_{m,w}(v_1,v')+k_{m,w}(v_1,v')k_{m,w}(v',v'').
\end{aligned}
\end{equation*}
The first two differences leads to a small contribution to $\mathcal{J}_4$

\begin{equation}\label{3.2.29}
  \frac{C}{N_0}\sup_{0\leq s\leq t}\|h_{R,\varepsilon}\|_\infty.
\end{equation}
For the main contribution of $k_{w,m}(v_1,v')k_{w,m}(v_1',v'')$. Using \eqref{3.2.3}, \eqref{3.2.27} reduces to
  \begin{equation}
    \begin{aligned}
        & \frac{1-\alpha}{\varepsilon^2} \mathbbm{1}_{t_1>0,t_1'>0} \int_{\substack{v_1\cdot n>0, \\ |v_1|\leq N_0}} \int_{0}^{t_1} \int_{0}^{t_1'-\varepsilon k_*} \int_{|v'|\leq 2N_0} \int_{|v''|\leq 2N_0} e^{\frac{\nu(s-t)}{\varepsilon}} e^{\frac{\nu(s'-t'_0)}{\varepsilon}} \\
        & \quad \times k_{w,m}(v_1,v') k_{w,m}(v_1',v'') \mathcal{R}(v,v_1) \sqrt{2\pi\mu(v'_0)\mu(v'_1)} |v_1'\cdot n| \\
        & \quad \times h_{R,\varepsilon}(s',x_1'+(s'-t_1')v_1',v'') \, ds' \, ds \, dv' \, dv_1' \, dv_1 \, dv'' \\
        & + \frac{\alpha}{\varepsilon^2} \mathbbm{1}_{t_1>0,t_1'>0} \int_{\substack{v_1\cdot n>0, \\ |v_1|\leq N_0}} \int_{\substack{v_1'\cdot n>0, \\ |v_1'|\leq N_0}} \int_{0}^{t_1} \int_{0}^{t_1'-\varepsilon k_*} \int_{|v'|\leq 2N_0} \int_{|v''|\leq 2N_0} e^{\frac{\nu(s-t)}{\varepsilon}} e^{\frac{\nu(s'-t'_0)}{\varepsilon}} \\
        & \quad \times k_{w,m}(v_1,v') k_{w,m}(v_1',v'') \mathcal{R}(v,v_1) \sqrt{2\pi\mu(v'_0)\mu(v'_1)} |v_1'\cdot n| \\
        & \quad \times h_{R,\varepsilon}(s',x_1'+(s'-t_1')v_1',v'') \, ds' \, ds \, dv' \, dv_1' \, dv_1 \, dv'' \\
        & := \mathcal{J}_{4,1} + \mathcal{J}_{4,2}.
    \end{aligned}
    \end{equation}
Notice that $v_0'=v'$. To estimate $\mathcal{J}_{4,1}$, we make a change of variable $v'\to y=x_1'+(s'-t_1')R_xv'.$ One has
  \begin{equation*}
    \mid \det\left(\frac{\partial y}{\partial v'}\right)\mid \geq (\varepsilon k_*)^3>0 \text{ for } s'\in[0,t_1'-k_*\varepsilon],
  \end{equation*}
  which yields that
  \begin{equation*}
  \begin{aligned}
      &\int_{|v'|\leq N_0}|h_{R,\varepsilon}(s',x_1'+(s'-t_1')R_xv',v'')|dv_1'\\
      \leq& C_{N_0}\left(\int_{|v'|\leq 2N_0}\mathbbm{1}_{\mathbb{R}^3_+}\left(x_1'+(s_1'-t_1')R_xv'\right)|h_{R,\varepsilon}(s',x_1'+(s'-t_1')R_xv',v'')|^2dv'\right)^{\frac{1}{2}}\\
      \leq& \frac{C_{N_0}}{(\varepsilon k_*)^{\frac{3}{2}}}\left(\int_{\mathbb{R}^3_+}|h_{R,\varepsilon}(s',y,v'')|^2dy\right)^\frac{1}{2}.
  \end{aligned}
  \end{equation*}
  Together with the definition of $f$ and $h_{R,\varepsilon}$, $\mathcal{J}_{4,1}$ can be bounded by $\frac{C_{N_0}}{(\varepsilon k_*)^{\frac{3}{2}}}\sup_{s\in [0,t]}\|f_{R,\varepsilon}\|_2$. Similarly, we make a change of variable $v'_1\to y_1=x_1'+(s'-t_1')v_1'$ and following the analogously arguments as above, $\mathcal{J}_{4,2}$ can be bounded by $\frac{C_{N_0}}{(\varepsilon k_*)^{\frac{3}{2}}}\sup_{s\in [0,t]}\|f_{R,\varepsilon}\|_2$.
One thereby has
\begin{equation}\label{3.2.30}
  \mathcal{J}_4\leq \frac{C}{N_0}\sup_{0\leq s\leq t}\|h_{R,\varepsilon}(s)\|_\infty+\frac{C_{N_0,k_*}}{\varepsilon^\frac{3}{2}}\sup_{0\leq s\leq t}\|f_{R,\varepsilon}\|_2.
\end{equation}
Collecting the above estimates, one has
\begin{equation}\label{3.2.31}
\begin{aligned}
   \sup_{0\leq s\leq t} \|h_{R,\varepsilon}(s)\|_\infty\leq &C\|h_{R,\varepsilon}(0)\|_\infty+C\sqrt{\varepsilon}^5 \sup_{0\leq s\leq t}\|h_{R,\varepsilon}(s)\|_\infty^2+C(\sqrt{\varepsilon}+\frac{1}{N_0}+\kappa_*)\sup_{0\leq s\leq t}\|h_{R,\varepsilon}(s)\|_\infty\\
  & +C\sqrt{\varepsilon}^3\sup_{0\leq s\leq t}\|h_{R,\varepsilon}(s)\|_\infty^2
    +\frac{C_{N_0,k_*}}{\varepsilon^\frac{3}{2}}\sup_{0\leq s\leq t}\|f_{R,\varepsilon}\|_2+C\sqrt{\varepsilon}^2.
\end{aligned}
\end{equation}
which concludes the Lemma \ref{Lm3.2} by choosing sufficiently small $\varepsilon>0, m > 0, \kappa_*> 0$ and
large $N_0 >0$.

\appendix
\section{Isotropic Functions}\label{App.A}
In this appendix, we list some significant functions that have been introduced in the boundary layer analysis.
The Burnett functions $A(v)$ and $B(v)$ are defined by
\begin{equation}\label{A.1}
  A_{ij}(v)=\left(v_iv_j-\frac{\delta_{ij}}{3}|v|^2\right)\sqrt{\mu},\quad B_i(v)=\left(\frac{|v|^2-5}{2}\right)v_i\sqrt{\mu}, \quad i,j=1,2,3.
\end{equation}
We also define the scalar function $C(v)$ by
\begin{equation}\label{A.2}
  C(v)=\left(\frac{1}{4}|v|^4-\frac{5}{2}|v|^2+\frac{15}{4}\right)\sqrt{\mu}.
\end{equation}
It is direct to see that $A_{i j},  {B}_{i}, C \in \mathcal{N}^{\perp}$ and $C$ is perpendicular to each entry of $A$ and $B$. We define

\begin{equation}\label{A.3}
  \begin{aligned}
  & \kappa_{1}:=\left\langle A _{i j}, \mathcal{L}^{-1} A_{i j}\right\rangle>0,\quad  i \neq j, \\
  & \kappa_{2}:=\left\langle B_{i}, \mathcal{L}^{-1} B_{i}\right\rangle>0.
  \end{aligned}
\end{equation}It's well-known that there exist two scalar functions $\alpha(|v|)$ and $\beta(|v|)$ such that (see \cite{desvillettes1994remark}, for instance)
\begin{equation}
  \hat{A}:=\alpha(|v|)A,\quad \hat{B}:=\beta(|v|)B,
\end{equation}
with
\begin{equation}
  \mathcal{L}\hat{A}=A,\quad \mathcal{L}\hat{B}=B.
\end{equation}
Based on \cite[Appendix A.2.9]{sone2007molecular} and \cite[Appendix B]{takata2012asymptotic}, we give below the definitions of the isotropic functions.
\begin{equation}
  \begin{aligned}
    &\mathcal{L} \left[(v_i\delta_{jk}+v_j\delta_{ik}+v_k\delta_{ij})D_1(|v|)+v_iv_jv_kD_2(|v|)\right] \\
    &=v_iv_jv_k\alpha(|v|)-\kappa_1(v_i\delta_{jk}+v_j\delta_{ik}+v_k\delta_{ij}),\\
    &\mathcal{L}(v_iv_3F_1(|v|))=v_i\hat{B}_3-\frac{1}{5}\kappa_2v_i,\\
    &\mathcal{L}F_2(|v|,v_3)=v_3\hat{B}_3-\kappa_2\frac{|v|^2-3}{3} .
  \end{aligned}
\end{equation}

\section{Axial Symmetry of Operators $\mathcal{L}$, $\mathcal{K}$ and $ {\Gamma}$ }\label{App.B}
For any function $f(v_3,|v|)$, the scattering and collision operator $\mathcal{K}$ and $\mathcal{L}$ satisfy the following axial symmetry:
\begin{equation}\label{B.1}
  J[v_if]=v_iJ_1[f], \quad (i=1,2),
\end{equation}
\begin{equation}\label{B.2}
  J[v_iv_jf]=v_iv_jJ_2[f]+J_3[f]\delta_{ij}, \quad (i,j=1,2 ),
\end{equation}
\begin{equation}\label{B.3}
  J[v_iv_jv_lf]=v_iv_jv_lJ_4+(v_i\delta_{jl}+v_j\delta_{il}+v_l\delta_{ij})J_5,\quad (i,j,l=1,2),
\end{equation}
where $J=\mathcal{L}$ and $\mathcal{K}$, and $J_1[f],J_2[f],J_3[f],J_4[f],J_5[f]$ are all functions of $v_3$ and $|v|$.

\begin{proof}[Sketch of proof:]
   One can get \eqref{B.1}-\eqref{B.3} for the operator $\mathcal{K}$ by a direct calculation. As for the linearized operator $\mathcal{L}$, equation \eqref{B.1}-\eqref{B.3} is a consequence of Appendix B.4.2 of \cite{sone2002kinetic} and Appendix A.2 of \cite{sone2007molecular}.
\end{proof}
Further, one can immediately obtain from \eqref{B.2} that:
\begin{equation}\label{B.4}
  J[\left(v_1^2-v_2^2\right)f]=\left(v_1^2-v_2^2\right)J_2[f].
\end{equation}
Note that we can divide any function $v_1^2f(v_3,|v|)$ into:
\begin{equation}\label{B.5}
  v_1^2f(v_3,|v|)=\frac{v_1^2-v_2^2}{2}f(v_3,|v|)+\frac{|v|^2-v_3^2}{2}f(v_3,|v|),
\end{equation}
where the latter term depends on $v_3,|v|$ only.

Similarly, from \eqref{B.3} we have
\begin{equation}\label{B.6}
  J\left[v_iv_jv_lf-\frac{1}{4}(\delta_{ij}+\delta_{jl}+\delta_{il})v_l(|v|^2-v_3^2)f\right]=\left[v_iv_jv_lf-\frac{1}{4}(\delta_{ij}+\delta_{jl}+\delta_{il})v_l(|v|^2-v_3^2)f\right]J_4[f].
\end{equation}
Similarly, for the bilinear term $\Gamma$, we have
\begin{equation}\label{B.7}
  \Gamma(v_if_1,f_2)=v_i\Gamma_1(f_1,f_2) \quad (i=1,2),
\end{equation}
and
\begin{equation}\label{B.8}
  \Gamma(v_if_1,v_jf_2)=v_iv_j\Gamma_2(f_1,f_2)+\delta_{ij}\Gamma_3(f_1,f_2), \quad (i,j=1,2)
\end{equation}
Formulas \eqref{B.4}-\eqref{B.8} will be used in the construction of $F^{bb}_k$ $(4\leq k\leq6)$ following the analogously arguments in subsection \ref{se1.4.3}. We omit the details here for simplicity.

\bibliography{reff} % 导入lib，ref为“ref.lib"的文件名
\bibliographystyle{abbrv}
\end{document}